\documentclass[11pt, lmargin=1.25in, rmargin=1.25in]{article}

\RequirePackage{amsthm,amsmath,amsfonts,amssymb}
\RequirePackage[authoryear]{natbib} 
\RequirePackage[colorlinks,citecolor=blue,urlcolor=blue,breaklinks]{hyperref}

\usepackage{bm}
\usepackage{mathrsfs}  
\usepackage{appendix}
\usepackage{commath}
\usepackage{dsfont}
\usepackage{subfig}
\usepackage{xcolor}
\usepackage{mathtools}
\usepackage{enumitem}
\usepackage{fancyhdr}
\usepackage{mathrsfs} 
\usepackage{geometry}
\usepackage{caption}
\usepackage{booktabs}
\usepackage{multirow}
\usepackage{lscape}

\usepackage{setspace}
\newcommand{\N}{\mathbb{N}}
\newcommand{\R}{\mathbb{R}}
\newcommand{\eps}{\varepsilon}
\renewcommand{\P}{\mathbb{P}}
\newcommand\E{\mathbb{E}}

\theoremstyle{plain}
\newtheorem{theorem}{Theorem}[section]
\newtheorem{lemma}[theorem]{Lemma}

\theoremstyle{remark}
\newtheorem{remark}{Remark}[section]

\newtheorem{assumption}{Assumption}[section]

\begin{document}

\title{Winsorized mean estimation with heavy tails and adversarial contamination\footnote{We thank two referees for helpful comments and suggestions.}}

\author{
\begin{tabular}{c}
Anders Bredahl Kock\footnote{
Kock's research was supported by the European Research Council (ERC) grant number 101124535 -- HIDI (UKRI EP/Z002222/1).  He is also a member of, and grateful for support from, i) the Aarhus Center for Econometrics (ACE), funded by the Danish National Research Foundation grant number DNRF186,  and ii) the Center for Research in Energy: Economics and Markets (CoRe).} \\ 
\footnotesize	University of Oxford \\
\footnotesize Department of Economics\\
\footnotesize	10 Manor Rd, Oxford OX1 3UQ
\\
\footnotesize	{\footnotesize	\href{mailto:anders.kock@economics.ox.ac.uk}{anders.kock@economics.ox.ac.uk}} 
\end{tabular}
\begin{tabular}{c}
David Preinerstorfer \\ 
{\footnotesize WU Vienna University of Economics and Business} \\
{\footnotesize Institute for Statistics and Mathematics} \\
{\footnotesize Welthandelsplatz 1, 1020 Vienna} \\ 
{\footnotesize	 \href{mailto:david.preinerstorfer@wu.ac.at}{david.preinerstorfer@wu.ac.at}}
\end{tabular}
}

\date{First version: April, 2025 \\
This version: June, 2026}

\maketitle	

\begin{abstract}
Finite-sample upper bounds on the estimation error of a winsorized mean estimator of the population mean in the presence of heavy tails and adversarial contamination are established. In comparison to existing results, the winsorized mean estimator we study avoids a sample-splitting device and winsorizes substantially fewer observations, which improves its applicability and practical performance.  
\end{abstract}

\section{Introduction}

Estimating the mean~$\mu$ of a distribution~$P$ on~$\R$ based on an i.i.d.~sample $X_1,\hdots,X_n$ is one of the most fundamental problems in statistics. It has long been understood that the sample average does not perform well in the presence of heavy tails or outliers. Sparked by the work of~\cite{catoni2012challenging}, recent years have witnessed much attention to the construction  of estimators~$\hat{\mu}_{n, \delta}=\hat{\mu}_{n, \delta}(X_1,\hdots,X_n)$ of~$\mu$ that exhibit finite-sample sub-Gaussian concentration even when~$P$ is heavy-tailed in the sense of possessing only two (finite) moments: that is, there exists an~$L\in(0,\infty)$, such that for all~$\delta\in(0,1)$ and~$n\in\N$
\begin{equation*}
	\envert[0]{\hat{\mu}_{n, \delta}-\mu}
	\leq
	L\sigma_2\sqrt{\frac{\log(2/\delta)}{n}}
\end{equation*}
with probability at least~$1-\delta$ and where $\sigma_2^2=E\del[1]{X_1-\mu}^2$.
The sample average does not exhibit such sub-Gaussian concentration, but other estimators (possibly depending on~$\delta$) have been constructed in, e.g.,~\cite{lerasle2011robust},~\cite{catoni2012challenging}, \cite{DLLO}, \cite{lugosi2019sub}, \cite{cherapanamjeri2019fast}, \cite{hopkins2020mean}, \cite{lee2022optimal}, \cite{minsker2023efficient}, \cite{gupta2024beyond}, \cite{gupta2024minimax}, \cite{minsker2024geometric}. Papers concerned with estimating the mean of a distribution on~$\R^d$ for~$d$ (much) larger than~$1$ often pay particular attention to constructing estimators that can be computed in (nearly) linear time. We refer to the overview in~\cite{lugosi2019mean} for further references and discussion on estimators with sub-Gaussian concentration properties.

Other works have studied estimators that are robust against \emph{adversarial contamination}: In this setting, an adversary inspects the sample~$X_1,\hdots,X_n$ and returns a corrupted (or contaminated) sample~$\tilde{X}_1,\hdots,\tilde{X}_n$ to the statistician, which estimators take as input. Thus, the \emph{identity} of the corrupted observations (or ``outliers")
\begin{equation*}
	\mathcal{O}=\mathcal{O}(X_1,\hdots,X_n):=\cbr[1]{i\in\cbr[0]{1,\hdots,n}:\tilde{X}_i\neq X_i},
\end{equation*}
as well as their \emph{values}, i.e.,~the values of~$\cbr[0]{\tilde{X}_i}_{i\in \mathcal{O}}$, can (but need not) depend on the uncontaminated~$X_1,\hdots,X_n$. In particular,~$\mathcal{O}$ can be a random subset of~$\cbr[0]{1,\hdots,n}$, and the adversary can use further external randomization in specifying~$\mathcal{O}$ and $\cbr[0]{\tilde{X}_i}_{i\in \mathcal{O}}$. We assume that at most~$\eta n$ of the contaminated observations~$\tilde{X}_1, \hdots, \tilde{X}_n$ differ from the uncontaminated ones, that is
\begin{equation}\label{eq:contamfrac}
	|\mathcal{O}(X_1,\hdots,X_n)|\leq \eta n,
\end{equation}
where~$\eta \in [0, 1]$ is non-random.\footnote{Note that (with the exception of the results on adaptation in Section~\ref{sec:adaptive})~$\eta$ need not be the smallest non-random number satisfying~\eqref{eq:contamfrac}.} The construction of estimators that are robust to adversarial contamination (and sometimes also heavy tails) along with finite-sample upper bounds on their error has been studied in, e.g., \cite{lai2016agnostic}, \cite{cheng2019high}, \cite{diakonikolas2019robust}, \cite{hopkins2020robust}, \cite{LM21}, \cite{minsker2021robust}, \cite{bhatt2022minimax}, \cite{depersin2022robust}, \cite{dalalyan2022all}, \cite{minasyan2023statistically}, \cite{minsker2023efficient}, \cite{oliveira2025finite}. The recent book by~\cite{diakonikolas2023algorithmic} provides further references and discussion of different contamination settings.

\cite{LM21} have shown that a sample-split based winsorized\footnote{\cite{LM21} refer to the estimator in Section 2 of their paper as a (modified) trimmed mean estimator, but it would perhaps be more common in the literature to call it a (modified) winsorized mean estimator and we hence do so.} mean estimator has sub-Gaussian concentration properties in an adversarial contamination setting.\footnote{We stress that the construction of estimators that make efficient use of the data in dimension one is not the main focus of~\cite{LM21}. Instead they focus on constructing estimators that depend optimally, in terms of rates, on the confidence level and the sample size in higher dimension.} The multivariate case was studied as well. In the present paper, we focus on the univariate case and use the ideas in \cite{LM21} to establish sub-Gaussian concentration properties under adversarial contamination for a winsorized mean estimator that removes some practical limitations of that analyzed in~\cite{LM21}:
\begin{itemize}
	\item The winsorized mean estimator we study does not require a sample split to determine the winsorization points. This allows for more efficient use of the data and makes the estimator permutation invariant.
	\item Whereas the estimator in~\cite{LM21} requires~$8\eta<1/2$, i.e.,~$\eta<1/16$, the estimator we analyze requires~$\eta<1/2$, thus extending the amount of contamination that is allowed. 
	\item The estimator we study only winsorizes slightly more than the smallest and largest~$\eta n$ observations, whereas the estimator analyzed in~\cite{LM21} winsorizes substantially more observations, which may be practically undesirable when it is known that at most~$\eta n$ observations have been contaminated.
\end{itemize} 

We provide upper bounds for any given number of moments~$m\in[1,\infty)$ that the uncontaminated observations possess. Typically, e.g., in \cite{LM21}, the focus is on the perhaps most important case~$m = 2$, but the flexibility in~$m$ is instrumental in~\cite{HDGauss}, where high-dimensional Gaussian and bootstrap approximations to the distribution of vectors of winsorized means under minimal moment conditions are established. In Section~\ref{sec:unknowneta} we study the setting where the statistician  knows an~$\eta$ that satisfies~\eqref{eq:contamfrac}. Since the smallest~$\eta$ for which~\eqref{eq:contamfrac} holds is typically unknown, Section~\ref{sec:adaptive} shows how a standard application of Lepski's method can be used to construct an estimator that adapts to that quantity. Section~\ref{sec:dependentdata} outlines the possibilities and challenges in extending our results to dependent data, and Section~\ref{sec:sims} contains numerical results comparing the winsorized mean to a range of other estimators.

\subsection{Data generating process}\label{sec:fin}
As outlined above, an adversary inspects the i.i.d.~sample~$X_1, \hdots, X_n$ from the distribution~$P$, corrupts at most~$\eta n$ of its values, and then gives the corrupted sample~$\tilde{X}_1,\hdots,\tilde{X}_n$ satisfying~\eqref{eq:contamfrac} to the statistician, who wants to estimate the mean of the (unknown) distribution~$P$. We summarize this, together with some assumptions, for later reference:

\begin{assumption}\label{ass:setting}
	The random variables~$X_1,\hdots,X_n$ are i.i.d.~with~$\E |X_1|^m<\infty$ for some~$m\in[1,\infty)$,~$\mu:=\E X_1$, and~$\sigma_m^m:=\E|X_1-\mu|^m$. The actually observed adversarially contaminated random variables are denoted by~$\tilde{X}_1,\hdots,\tilde{X}_{n}$ and satisfy~\eqref{eq:contamfrac}.	
\end{assumption}
Note that although the uncontaminated data~$X_1,\hdots,X_n$ is assumed to be independent, this need not to be the case for the contaminated data~$\tilde{X}_1,\hdots,\tilde{X}_n$.

\section{Performance guarantees for known~$\eta$}\label{sec:unknowneta} 
We first study winsorized mean estimators requiring knowledge of~$\eta$. To this end, for real numbers~$x_1,\hdots,x_n$, we denote by~$x_1^*\leq \hdots\leq x_n^*$ their non-decreasing rearrangement. Let~$-\infty<\alpha\leq\beta<\infty$ and define
\begin{equation}\label{eq:winsor}
	\phi_{\alpha,\beta}(x)
	:=
	\begin{cases}
		\alpha\qquad \text{if }x<\alpha\\
		x\qquad \text{if }x\in[\alpha,\beta]\\
		\beta\qquad \text{if }x>\beta.
	\end{cases}
\end{equation} 
For~$\eps\in(0,1/2]$, let~$\hat{\alpha}=\tilde X_{\lceil \eps n \rceil}^*$ and $\hat{\beta}=\tilde X_{\lfloor(1-\eps )n\rfloor+1}^*$.\footnote{We consider~$\eps\in(0,1/2]$ since otherwise~$\hat{\alpha}$ could exceed~$\hat{\beta}$. Note that~$\hat{\mu}_n$ is a sample median for~$\varepsilon = 1/2$.} We consider winsorized estimators of the form
\begin{equation}\label{eqn:winsmean}
	\hat{\mu}_n=\hat{\mu}_n(\eps):=\frac{1}{n}\sum_{i=1}^n\phi_{\hat\alpha,\hat\beta}(\tilde{X}_i),
\end{equation}
Under adversarial contamination it is clear that any such estimator can perform arbitrarily badly unless at least the smallest and largest~$\eta n$ observations are winsorized. Thus, one must choose~$\eps\geq \eta$, implying in particular that~$\eta\leq1/2$ must hold.\footnote{Any estimator breaks down if half of the sample (or more) is (adversarially) contaminated, so it is no real restriction to focus on the case where~$\eta < 1/2$.} For a desired ``confidence level''~$\delta\in(0,1)$, we choose~$\eps$ as
\begin{equation}\label{eq:epsfam}
	\eps=
	\eps(\eta)
	:=
	\lambda_1\cdot \eta +\lambda_2\cdot \frac{\log(6/\delta)}{n}
	,\qquad \text{for fixed }\lambda_1\in(1,\infty)\text{ and }\lambda_2\in (0,\infty).
\end{equation} 

The estimator~$\hat{\mu}_n$ resulting from this choice of~$\eps$ is similar to the winsorized mean estimator in~\cite{LM21}. However, their approach uses a sample split to calculate~$\hat{\alpha}$ and $\hat{\beta}$ on one half of the sample and then computes the average in~\eqref{eqn:winsmean} only over the other half. This has the effect of ``halving'' the sample size and leads to an estimator that is not permutation invariant. Furthermore, their estimator corresponds to choosing~$\lambda_1=8$ and (essentially)~$\lambda_2=24$ above (note that their~$N$ is our~$n/2$ due to their sample split). As a consequence, their~$\eps$ exceeds~$1/2$ for many values of~$(n,\eta,\delta)$, rendering their estimator unimplementable, cf.~Section~\ref{sec:sims}. Furthermore, whenever their~$\eps\in(0,1/2]$, this implies that~$\eta<\eps/8\leq 1/16$, such that at most~$6.25\%$ of the observations can be adversarially contaminated in their implementation. It may be inefficient use of the data to use a sample split, and to winsorize (slightly more than) the smallest and largest $8\eta$ fraction of the remaining observations if one knows that at most~$\eta n$ observations are contaminated. Our implementation only winsorizes (slightly more than) the~$\lambda_1\eta n$ smallest and largest observations, and we recommend choosing~$\lambda_1$ only slightly larger than~$1$, e.g.,~$\lambda_1=1.01$.  Concerning the choice of~$\lambda_2$, the simulations in Section~\ref{sec:sims} suggest that small values of~$\lambda_2$ such as~$\lambda_2=0.2$ work well.

Our theoretical guarantees below for~$\hat{\mu}_n(\eps)$ apply for any $\eps$ in \eqref{eq:epsfam} satisfying 	
\begin{equation}\label{eq:epscond}
	2\eps +\frac{\log(6/\delta)}{n}+\sqrt{\del[2]{\frac{\log(6/\delta)}{n}}^2+4\frac{\log(6/\delta)}{n}\eps}<1.
\end{equation}
Note that this condition implies~$\eta < \eps < 1/2$. Although~\eqref{eq:epscond} is stronger than imposing~$\eps\in(0,1/2]$, which is all that is needed to \emph{implement}~$\hat{\mu}_n$ in~\eqref{eqn:winsmean}, note that~$\log(6/\delta)/n$ in~\eqref{eq:epscond} is typically small. Thus, for large~$n$ the requirement on~$\eps$ in~\eqref{eq:epscond} essentially reduces to~$\eps\in(0,1/2)$. In the special case of~$\eta=0$, such that~$\eps=\lambda_2\cdot \log(6/\delta)/n$,~\eqref{eq:epscond} reduces to
\begin{equation*}
	\del[1]{2\lambda_2+1+\sqrt{1+4\lambda_2}}\frac{\log(6/\delta)}{n}<1,
\end{equation*}
which is typically satisfied (even for moderate~$n$) if~$\lambda_2$ is small.
\begin{remark}
	Actually, the condition in~\eqref{eq:epscond} is just a conservative (simple) sufficient condition for the following milder condition that one could also work with (we have chosen not to, because it is more cumbersome and difficult to interpret): Writing~$$A_+=1-\lambda_1^{-1}\mathds{1}\cbr[0]{\eta>0}\in(0,1] \quad \text{ and } \quad A_-=1+\lambda_1^{-1}\mathds{1}\cbr[0]{\eta>0}\in [1,\infty),$$ and denoting by~$W_0$ and~$W_{-1}$ the principal and lower branch of Lambert's~$W$ function (cf., e.g., \cite{Corless1996}), respectively,~\eqref{eq:epscond} could be replaced by
	\begin{equation}\label{eqn:altcond}
		\eps\del[3]{-A_+W_0\del[1]{-e^{-(\frac{\log(6/\delta)}{\eps n}+A_+)/A_+}}-A_-W_{-1}\del[1]{-e^{-(\frac{\log(6/\delta)}{\eps n}+A_-)/A_-}}}<1.
	\end{equation}
	By \eqref{eq:welldefined} of Lemma~\ref{lem:cControl} in the appendix, the left-hand side of~\eqref{eqn:altcond} is upper bounded by the left-hand side of~\eqref{eq:epscond}, leading to the condition in~\eqref{eq:epscond}. Note, however, that the latter condition implies~$\log(6/\delta)/n < 1$, which is repeatedly used in the proofs.
\end{remark}

We next present an upper bound on the estimation error of~$\hat{\mu}_{n}(\eps(\eta))$ as defined in~\eqref{eqn:winsmean}; note that the notation emphasizes the dependence of the estimator on~$\eta$ to set it apart from the estimator adapting to the smallest~$\eta$ satisfying~\eqref{eq:contamfrac} studied in Section~\ref{sec:adaptive}. 
\begin{theorem}\label{thm:maintext}
	Fix~$n\in\N$,~$\delta\in(0,1)$, and let Assumption~\ref{ass:setting} be satisfied with~$m\in[1,\infty)$. Let~$\lambda_1\in(1,\infty)$ and~$\lambda_2\in(0,\infty)$. There exist positive constants~$\mathfrak{A}_m(\lambda_1, \lambda_2)$ and~$\mathfrak{B}_m(\lambda_1, \lambda_2)$, depending only on~$\lambda_1$,~$\lambda_2$, and~$m$, such that if~$\eps(\eta)$ is chosen as in~\eqref{eq:epsfam} and satisfies~\eqref{eq:epscond}, then, with probability at least~$1-\delta$, we have
	\begin{equation}\label{eqn:mbound}
		\envert[1]{\hat{\mu}_{n}(\eps(\eta))-\mu}
		\leq 
		\sigma_m\del[4]{\mathfrak{A}_m(\lambda_1, \lambda_2) \cdot \eta^{1-\frac{1}{m}}+\mathfrak{B}_m(\lambda_1, \lambda_2) \cdot \del[2]{\frac{\log(6/\delta)}{n}}^{1-\frac{1}{m\wedge 2}}},
	\end{equation}
	which, in case~$m = 2$, simplifies to
	\begin{equation}\label{eqn:m2bound}
		\envert[1]{\hat{\mu}_{n}(\eps(\eta))-\mu}
		\leq 
		\sigma_2\del[3]{\mathfrak{A}_2(\lambda_1, \lambda_2) \cdot  \sqrt{\eta}+\mathfrak{B}_2(\lambda_1, \lambda_2) \cdot \sqrt{\frac{\log(6/\delta)}{n}}}.
	\end{equation}
	[The constants~$\mathfrak{A}_m(\lambda_1, \lambda_2)$ and~$\mathfrak{B}_m(\lambda_1, \lambda_2)$ are explicitly given in the proof.]\footnote{In case~$m = 1$ and~$\eta = 0$ one can set~$\eta^{1-1/m} = 0$ in the upper bound.\label{foot:m1eta0}}
\end{theorem}
The dependence of~\eqref{eqn:mbound} on~$\eta$ appears to be optimal up to multiplicative constants for all~$m\in[1,\infty)$. This follows from the argument on pages 396--397 in~\cite{LM21} upon replacing $\sqrt{\eta}$ by~$\eta^{1-1/m}$ and~$\sigma_X$ by~$\sigma_m$, respectively, in the distribution constructed in the remark on their page~397. 

Larger~$m$ correspond to lighter tails of the~$X_1,\hdots,X_n$. This makes it easier to classify large contaminations as outliers, which, essentially, ``restricts'' the meaningful contamination strategies of the adversary. Thus, it is sensible that larger~$m$ lead to a better dependence on the contamination rate~$\eta$. 

The proof of Theorem~\ref{thm:maintext} builds on a decomposition of the estimation error outlined in Appendix \ref{sec:Decomposition}. A similar decomposition was implicitly used in~\cite{LM21}. However, in contrast to~\cite{LM21}, we do not use a sample split to determine the winsorization locations~$\hat{\alpha}=\tilde X_{\lceil \eps n \rceil}^*$  and $\hat{\beta}=\tilde X_{\lfloor(1-\eps )n\rfloor+1}^*$. Furthermore, to reduce excessive winsorization, i.e., to allow~$\lambda_1\in(1,\infty)$ and~$\lambda_2\in(0,\infty)$ instead of~$\lambda_1=8$ and~$\lambda_2=24$ in \cite{LM21}, we carefully bound~$\hat{\alpha}$ and~$\hat{\beta}$ in Lemma~\ref{lem:quantiles}. These bounds are fundamental to our approach. We here exploit exponential concentration inequalities  tailored to the Binomial distribution (in particular the inequalities in Lemma~\ref{lem:Chernoff}, which are taken from \cite{hagerup1990guided}) rather than using the more ``general purpose'' Bernstein inequality (which the argument in \cite{LM21} is based on). To establish the feasibility of our approach, we first carefully study the exponent in these concentration inequalities and solutions to equations related to these that can be expressed in terms of Lambert's~$W$ function, cf.~Lemmas~\ref{lem:Prophfs} and~\ref{lem:cControl}. [We also note that if one replaces Lemmas~\ref{lem:Chernoff}--\ref{lem:cControl} by the Bernstein inequality and an analogous careful analysis of the corresponding exponent, this would result in the restriction~$\lambda_2\geq 2/3$ when~$\eta=0$, so that it is not possible to allow~$\lambda_2$ to take any value in~$(0,\infty)$ with that approach.]

\section{Adapting to the smallest~$\eta$ by Lepski's method}\label{sec:adaptive}

In practice, an~$\eta$ for which~\eqref{eq:contamfrac} holds is often unknown. Furthermore, even if one happens to know some~$\eta$ satisfying~\eqref{eq:contamfrac}, the upper bound established in Theorem~\ref{thm:maintext} increases (for~$m > 1$) in~$\eta$, so that one would like to choose~$\eta$ as small as possible. We now construct an estimator that adapts to the smallest non-random~$\eta$ for which~\eqref{eq:contamfrac} is satisfied, i.e., to
\begin{equation}\label{eq:etamin}
	\eta_{\min}:=\min\cbr[0]{\eta \in [0, 1] :|\mathcal{O}(X_1,\hdots,X_n)|/n\leq \eta}.
\end{equation}
The construction of this adaptive estimator is based on (the ideas underlying) Lepski's method, cf., e.g.,~\cite{lepskii1991problem, lepskii1992asymptotically, lepskii1993asymptotically}. Our specific implementation combines elements of the proofs of Theorem 3 in~\cite{dalalyan2022all} and Theorem 4.2 in~\cite{DLLO}. 

Fix~$m\in[1,\infty)$ as in Assumption~\ref{ass:setting}. In addition, let~$\rho\in(0,1)$ and suppose that~$\eta_{\min} \in [0,0.5\rho]$. For~$\delta>6\exp(-n/2)$ we define~$g_{\max}:=\lceil \log_\rho(2\log(6/\delta)/n)\rceil $ and the geometric grid of points~$\eta_j:=0.5\rho^j$ for~$j\in[g_{\max}]:=\cbr[1]{1,\hdots,g_{\max}}$. Let~$$g^*:=\max\cbr[0]{j\in [g_{\max}]: \eta_{\min} \leq \eta_j}.$$ Thus,~$\eta_{g^*}$ is the smallest~$\eta_j$ exceeding (the unknown)~$\eta_{\min}$. For~$x\in\R$ and~$r\in(0,\infty)$, let~$\mathbb{B}(x,r):=\cbr[0]{y\in\R:|y-x|\leq r}$. Furthermore, define for every~$z\in[0,\infty)$ the quantity (cf.~Theorem~\ref{thm:maintext} and its proof for explicit expressions for $\mathfrak{A}_m(\lambda_1, \lambda_2)$ and~$\mathfrak{B}_m(\lambda_1, \lambda_2)$)
\begin{equation*}
	B(z):=	\sigma_m \cdot \left(
	\mathfrak{A}_m(\lambda_1, \lambda_2) \cdot z^{1-\frac{1}{m}}+\mathfrak{B}_m(\lambda_1, \lambda_2) \cdot \del[2]{\frac{\log(6g_{\max}/\delta)}{n}}^{1-\frac{1}{m\wedge 2}}\right),
\end{equation*}
%
where, for notational convenience, we do not highlight the dependence of~$B$ on~$\sigma_m, \lambda_1, \lambda_2$ and~$m$. Recall that~$\delta\in(0,1)$, and let
\begin{equation}\label{eqn:epsAdef}
	\eps_A(\eta):=\lambda_1\cdot \eta +\lambda_2\cdot \frac{\log(6g_{\max}/\delta)}{n},\quad \text{for fixed }\lambda_1\in(1,\infty)\text{ and }\lambda_2\in (0,\infty);
\end{equation}
noting that~$\eps_A(\eta)$ corresponds to~$\eps(\eta)$ in~\eqref{eq:epsfam} with~$\delta$ there replaced by~$\delta/g_{\max}$. Define the analogue
\begin{equation}\label{eq:epscondLepski}
	2\eps_A(\eta) +\frac{\log(6g_{\max}/\delta)}{n}+\sqrt{\del[2]{\frac{\log(6g_{\max}/\delta)}{n}}^2+4\frac{\log(6g_{\max}/\delta)}{n}\eps_A(\eta)}<1
\end{equation}
to~\eqref{eq:epscond}; the difference (again) being that~$\delta$ in~\eqref{eq:epscond} is replaced by~$\delta/g_{\max}$ in~\eqref{eq:epscondLepski}. Finally, set
\begin{align*}
	\mathbb{I}(\eta_j)
	:=
	\begin{cases}
		\mathbb{B}\del[1]{\hat{\mu}_{n}(\eps_A(\eta_j)),B (\eta_j)}\qquad &\text{if } \eps_A(\eta_j)\text{ satisfies~\eqref{eq:epscondLepski}}\\
		\R	&\text{if } \eps_A(\eta_j)\text{ does not satisfy~\eqref{eq:epscondLepski}},
	\end{cases}
\end{align*}
for $j\in[g_{\max}]$, and define
\begin{equation*}
	\hat{g}:=\max\cbr[3]{g\in [g_{\max}]:\bigcap_{j=1}^g \mathbb{I}(\eta_j)\neq \emptyset}.
\end{equation*}
Under the assumptions of Theorem~\ref{thm:adapt}, $\bigcap_{j=1}^{\hat{g}} \mathbb{I}(\eta_j)$ will be shown to be a non-empty finite interval (possibly degenerated to a single point). Thus, we can define the estimator~$\hat{\mu}_{n,A}$ as the (measurable) midpoint of~$\bigcap_{j=1}^{\hat{g}} \mathbb{I}(\eta_j)$. Note that~$\hat{\mu} _{n,A}$ can be implemented \emph{without} knowledge of~$\eta_{\min}$. In addition,~$\hat{\mu}_{n,A}$ adapts to the unknown~$\eta_{\min}$ in the following sense.
\begin{theorem}\label{thm:adapt}
	Fix~$n\geq 4$,~$\delta\in (6\exp(-n/2),1)$, and let Assumption~\ref{ass:setting} be satisfied with~$m\in [1,\infty)$. Let~$\lambda_1\in(1,\infty)$ and~$\lambda_2\in(0,\infty)$. Furthermore, let~$\rho\in(0,1)$, suppose that $\eta_{\min}\in [0,0.5\rho]$, and that~$\eps_A(\eta_{g^*})$ as defined in~\eqref{eqn:epsAdef} satisfies~\eqref{eq:epscondLepski}. Let~$\mathfrak{A}_m(\lambda_1, \lambda_2)$ and~$\mathfrak{B}_m(\lambda_1, \lambda_2)$ be as in Theorem~\ref{thm:maintext} (cf.~also its proof), and set~$\mathfrak{C}_m(\lambda_1, \lambda_2) := \mathfrak{A}_m(\lambda_1, \lambda_2) + \mathfrak{B}_m(\lambda_1, \lambda_2)$. Then, with probability at least~$1-\delta$, we have (suppressing the dependence of~$\mathfrak{A}_m(\lambda_1, \lambda_2)$ and $\mathfrak{C}_m(\lambda_1, \lambda_2)$ on~$\lambda_1$ and~$\lambda_2$)\footnote{We adopt the same convention as in in Footnote \ref{foot:m1eta0}: In case~$m = 1$ and~$\eta_{\text{min}} = 0$ one can set~$\eta_{\text{min}}^{1-1/m} = 0$ in the upper bound in~\eqref{eq:UBadapt}.}
	\begin{eqnarray}\label{eq:UBadapt}
		\envert[1]{\hat{\mu}_{n,A}-\mu}
		\leq 
		2\sigma_m \cdot \left(
		\mathfrak{A}_m \cdot \left[\frac{\eta_{\min}}{\rho}\right]^{1-\frac{1}{m}}+\mathfrak{C}_m \cdot \del[2]{\frac{\log(6g_{\max}/\delta)}{n}}^{1-\frac{1}{m\wedge 2}}\right),
	\end{eqnarray}
	which, in case~$m = 2$, simplifies to
	\begin{equation*}
		\envert[1]{\hat{\mu}_{n,A}-\mu}
		\leq 
		2\sigma_2 \cdot \left(
		\mathfrak{A}_2 \cdot \sqrt{\frac{\eta_{\min}}{\rho} }+\mathfrak{C}_2 \cdot \sqrt{\frac{\log(6g_{\max}/\delta)}{n}}\right).
	\end{equation*}
\end{theorem}
The estimator~$\hat{\mu}_{n,A}$, which does \emph{not} have access to~$\eta_{\min}$, has the same dependence on~$\eta_{\min}$ (up to multiplicative constants) as the estimator~$\hat{\mu}_{n}(\eps(\eta_{\min}))$ from Theorem~\ref{thm:maintext} that \emph{knows}~$\eta_{\min}$ and uses~$\eta=\eta_{\min}$. However, observe that~$\hat{\mu}_{n,A}$ only adapts to~$\eta_{\min}\in[0,0.5\rho]\subsetneq [0,0.5)$. This gap in the adaptation zone can be made arbitrarily small by choosing~$\rho$ close to (but strictly less than) one. We also note that the terms in the upper bound in~\eqref{eq:UBadapt} that do not involve the fraction of contaminated observations are \emph{larger} than the corresponding terms in the upper bound in~\eqref{eqn:mbound}. This suggests that the adaptivity property of~$\hat{\mu}_{n,A}$ does not come ``for free'' and that one should not use the adaptive estimator if one (roughly) knows~$\eta_{\min}$.

We emphasize that~$\hat{\mu}_{n,A}$ incorporates knowledge of~$\sigma_m$. This can be avoided by replacing~$\sigma_m$ in the construction of~$\hat{\mu}_{n,A}$ (i.e., in the definition of~$B$) by an upper bound on it. The argument used to prove Theorem~\ref{thm:adapt} still goes through (with slight modifications) for this modified estimator, and establishes a similar statement as in~\eqref{eq:UBadapt}, but where~$\sigma_m$ has to be replaced by its upper bound.\footnote{It is common that an upper bound on~$\sigma_m$ or related quantities is needed when constructing estimators adapting to various quantities (such as~$\eta_{\min}$), cf., e.g., \cite{DLLO} or \cite{dalalyan2022all}.}

\begin{remark}
	The proof of Theorem~\ref{thm:adapt} shows that with probability at least~$1-\delta$ it holds that~$\hat{\mu}_{n,A}$ is within a distance~$B(\eta_{g^*})$ to the \emph{infeasible} estimator $\hat{\mu}_n(\varepsilon_A(\eta_{g^*}))$ that uses the \emph{unknown} smallest upper bound~$\eta_{g^*}$ on~$\eta_{\min}$ from the grid~$\cbr[0]{\eta_j:j\in[g_{\max}]}$. Thus, the adaptive estimator~$\hat{\mu}_{n,A}$ essentially works by selecting among the estimators $$\cbr[1]{\hat{\mu}_n(\varepsilon_A(\eta_j)): j\in[g_{\max}]}$$ from Theorem~\ref{thm:maintext} the one that uses the lowest value~$\eta_j$ exceeding~$\eta_{\min}$.
\end{remark}

\begin{remark}\label{rem:otherestim}
	At the price of larger multiplicative constants in the upper bound only, one could have defined the adaptive estimator as~$\tilde{\mu}_{n}=\hat{\mu}_n(\varepsilon_A(\eta_{\hat{g}}))$, which is an element of the grid of estimators~$\cbr[1]{\hat{\mu}_n(\varepsilon_A(\eta_j)): j\in[g_{\max}]}$, and thus arguably more natural than~$\hat{\mu}_{n,A}$. In Remark~\ref{rem:altadapt} in the appendix we establish an upper bound on $|\tilde{\mu}_{n}-\mu|$ similar to that in Theorem~\ref{thm:adapt}.
\end{remark}

\section{Dependent data}\label{sec:dependentdata}
In this section, we discuss the possibilities for --- and challenges involved in --- extending Theorem~\ref{thm:maintext} to dependent data.  Inspection of the proof of Theorem~\ref{thm:maintext} and the supporting lemmas leading to it shows that the dependence notion entertained should be ``stable'' under transformations applied to the individual observations such as winsorization and taking certain indicators. Furthermore, in the current method of proof, the independence of~$X_1,\hdots,X_n$ is used in establishing
\begin{enumerate}
	\item Lemma~\ref{lem:contdisc} to avoid imposing continuity of the cdf of the~$X_i$. 
	\item Lemma~\ref{lem:quantiles}, which provides control of the winsorization locations~$\hat{\alpha}=\tilde X_{\lceil \eps n \rceil}^*$  and $\hat{\beta}=\tilde X_{\lfloor(1-\eps )n\rfloor+1}^*$. Here we make use of Chernoff-bound based concentration inequalities tailored to the binomially distributed~$S_n=\sum_{i=1}^n \mathds{1}\del[1]{X_i\leq Q_{p}(X_1)}$ and related sums (Lemma \ref{lem:Chernoff}); for~$Q_p(X_1)=\inf\cbr[1]{z\in \R:\P(X_1\leq z)\geq p}$ for~$p\in(0,1)$. The feasibility of this approach relies on an analysis of the existence, uniqueness, and properties of solutions to equations related to the exponent of the Chernoff-bound in Lemmas~\ref{lem:Prophfs} and~\ref{lem:cControl}.
	\item Lemma~\ref{lem:trimmeanconc} via Bernstein's inequality for sums of  independent bounded random variables.	 
\end{enumerate}
A version of Lemma~\ref{lem:contdisc} can likely be established for some typical dependence concepts. Alternatively, one could also impose the~$X_i$ to have a continuous cdf (which, however, would limit the scope of the results). For these reasons, the first item does not constitute a major obstacle.

Since~$S_n$ defined in the second item of the above enumeration is a sum of bounded random variables, one can, in principle, replace the use of the Chernoff-bound for the Binomial distribution in Lemma~\ref{lem:quantiles} and the use of Bernstein's inequality in Lemma~\ref{lem:trimmeanconc} by a Bernstein inequality valid for the form of dependence that one is willing to entertain. For example,~\cite{merlevede2009bernstein, merlevede2011bernstein} have established Bernstein inequalities under geometric~$\alpha$-mixing, and, more recently,~\cite{Hang17} have established a Bernstein inequality for stochastic processes that include~$\phi$-mixing processes. Note, however, that already in the i.i.d.~case using only the Bernstein inequality leads to the unnecessary restriction~$\lambda_2\geq 2/3$ when~$\eta=0$, cf.~the discussion at the end of Section~\ref{sec:unknowneta}. This would carry over to the dependent case. 

Note also that Bernstein inequalities for dependent data often contain unknown population quantities such as mixing coefficients and ``long-run'' variances; the latter themselves being functions of unknown covariances, cf.~Theorems 1 and 2 in \cite{merlevede2009bernstein} and Theorem 1 in~\cite{merlevede2011bernstein}. Thus, to establish an analogue of Lemma~\ref{lem:Prophfs},~$\lambda_1$ and~$\lambda_2$ would likely have to be restricted in a way depending on these unknown quantities, making the practical implementation of the associated winsorized mean difficult. In addition, Bernstein inequalities for dependent data can involve (powers of) logarithmic terms not present in the Bernstein inequality for independent data, implying that the second summand in the definition of~$\eps$ in~\eqref{eq:epsfam} would possibly have to be chosen in a different  manner specific to the dependence notion employed. 

Hence, while our general approach can likely be extended also to dependent observations, the domains of~$\lambda_1$ and~$\lambda_2$ (as well as the specific form of~$\varepsilon$) will possibly have to be restricted, the restriction incorporating the dependence concept entertained. The resulting estimators could be of limited practical value, if they have to be based on large values for~$\lambda_1$ and~$\lambda_2$. We therefore leave a careful study of the dependent case to future research.

\section{Numerical evidence}\label{sec:sims}
In this section, we numerically investigate the performance of the winsorized mean estimators studied. Throughout, the winsorized mean~$\hat{\mu}_n$ in~\eqref{eqn:winsmean} with~$\eps(\eta)$ chosen as in~\eqref{eq:epsfam} is implemented with~$\lambda_1=1.01$ to avoid excessive winsorization. The sensitivity to the choice of~$\lambda_2$ is studied by implementing~$\hat{\mu}_n$ with~$\lambda_2\in\cbr[0]{0.2,0.5,1}$. 

The adaptive estimator~$\hat{\mu}_{n,A}$ from Section~\ref{sec:adaptive} is primarily a theoretical construction used to demonstrate that adaptation to the unknown~$\eta_{\min}$ is possible. Recall also that implementation of~$\hat{\mu}_{n,A}$ requires knowledge of~$m$ and~$\sigma_m$. With these caveats in mind, we implement~$\hat{\mu}_{n,A}$ with~$\lambda_1=1.5$ and~$\lambda_2=0.2$.\footnote{The constants~$\mathfrak{A}(\lambda_1, \lambda_2)$ and~$\mathfrak{B}(\lambda_1, \lambda_2)$  entering the definition of~$B(z)$ become very large as~$\lambda_1$ approaches one. We hence use~$\lambda_1=1.5$ and reiterate that the results for this estimator are illustrative only. $\lambda_2=0.2$ is used since this turns out to work quite well for~$\hat{\mu}_n$ on which~$\hat{\mu}_{n,A}$ is based.} For comparison, we also implement the sample average, the trimmed mean as in Theorem 1.3.1 in \cite{oliveira2025finite}, the winsorized mean from Section 2 in~\cite{LM21}, and the median-of-means estimator as in Theorem 2 in \cite{lugosi2019mean} (the latter being built for a setting that does not take into account adversarial contamination).

To assess the performance of winsorized and trimmed mean estimators it is useful to consider distributions for which the mean and median (here defined as the smallest~$1/2$-quantile of the cdf of~$X_1$) differ: Otherwise, estimators that winsorize or trim excessively and hence ``approach'' the empirical median (which concentrates strongly around the population median \emph{irrespective} of the number of moments the~$X_i$ possess, cf.~Lemma~\ref{lem:quantiles} in the appendix) may perform artificially well simply because the population median equals the population mean. To construct a simple example of such a distribution, denote by~$\delta_a$ the Dirac measure at~$a\in\R$ and by~$\mathsf{P}_{t,\gamma}$ the Pareto distribution with scale parameter~$t>0$ and shape~$\gamma>1$. The uncontaminated~$X_i$ are generated from the (mean-zero) mixture
\begin{equation*}
	\mathsf{m}=\mathsf{m}_{t,\gamma}=0.5\cdot\delta_{-b}+0.5\cdot \mathsf{P}_{t,\gamma}\ast \delta_{-b},
\end{equation*}
where
\begin{equation*}
b=b_{t,\gamma}=0.5\int x\mathsf{P}_{t,\gamma}(dx)=\frac{\gamma t}{2(\gamma-1)},
\end{equation*}
and~$\mathsf{P}_{t,\gamma}\ast \delta_{-b}$ is the convolution of~$\mathsf{P}_{t,\gamma}$ and~$\delta_{-b}$. Note that
\begin{enumerate}
	\item $\mathsf{m}$ possesses all moments strictly less than~$\gamma$, since the Pareto distribution~$\mathsf{P}_{t,\gamma}$ possesses all moments strictly less than~$\gamma$.
	\item the median of~$\mathsf{m}$ is~$-b=\frac{-\gamma t}{2(\gamma-1)}$, whereas the mean is~$0$. Thus, for any given number of moments that~$\mathsf{m}$ possesses (controlled via~$\gamma>1$), one can control the distance~$b$ between the mean and median via~$t$. 
\end{enumerate}
Throughout we use~$t=2$ and~$\gamma=m+0.01$ for~$m\in\cbr[0]{2,3}$ such that~$\mathsf{m}$ has only slightly more than~$m$ moments. All estimators use~$\delta=0.01$ and all simulations are based on~$100{,}000$ replications.  We consider~$n\in\cbr[0]{200,500}$. For the sake of comparison to~$X_i\sim \mathsf{m}_{t,\gamma}$, we also report some findings from simulations wherein~$X_i\sim \mathsf{t}(\gamma)$, the~$t$-distribution with~$\gamma$ degrees of freedom for $\gamma=m+0.01$ for~$m\in\cbr[0]{2,3}$. For these distributions the median equals the mean.

\subsection{No contamination:~$\eta_{\min}=0$} 
\begin{landscape}
\begin{figure}
	\includegraphics[width=10cm, height=7cm]{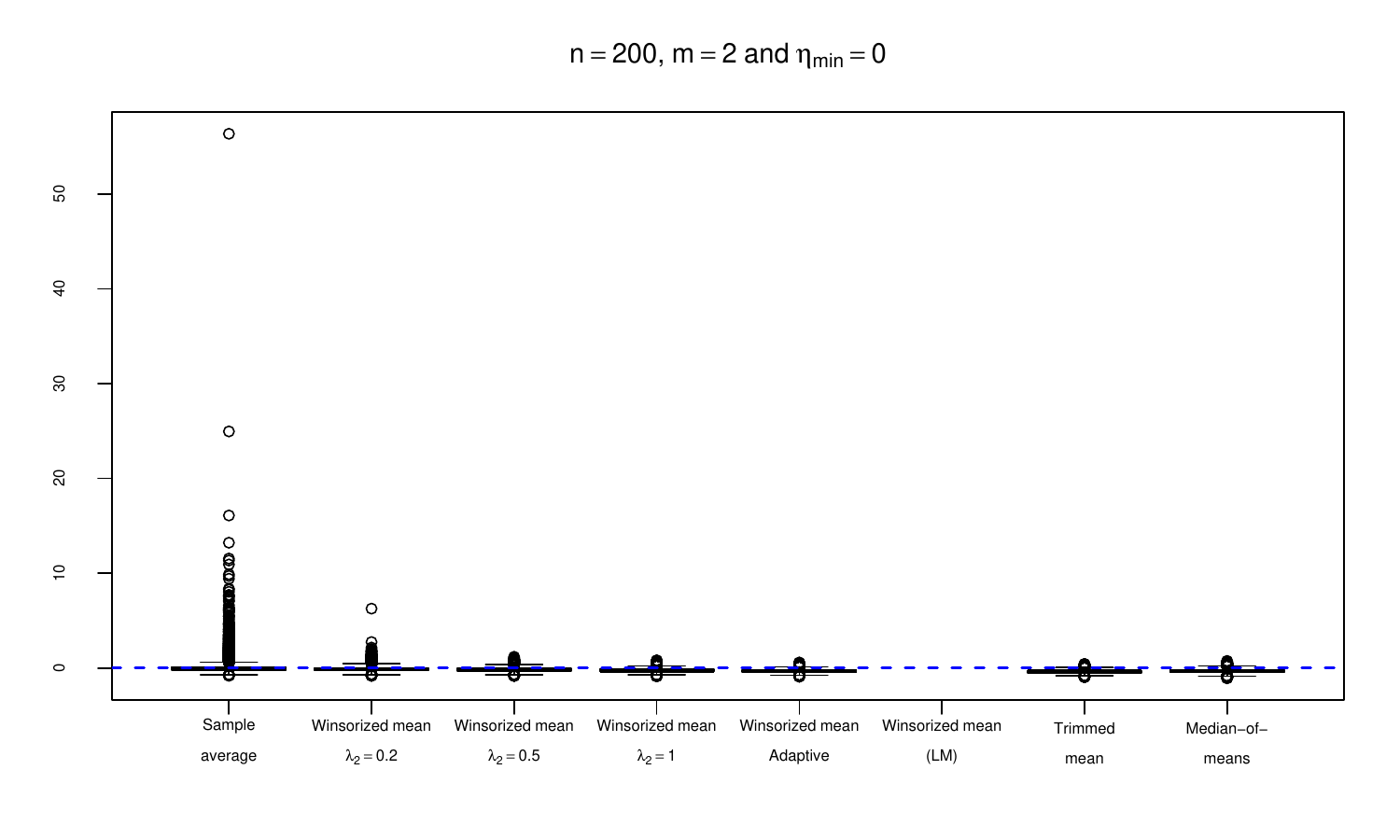}
	\includegraphics[width=10cm, height=7cm]{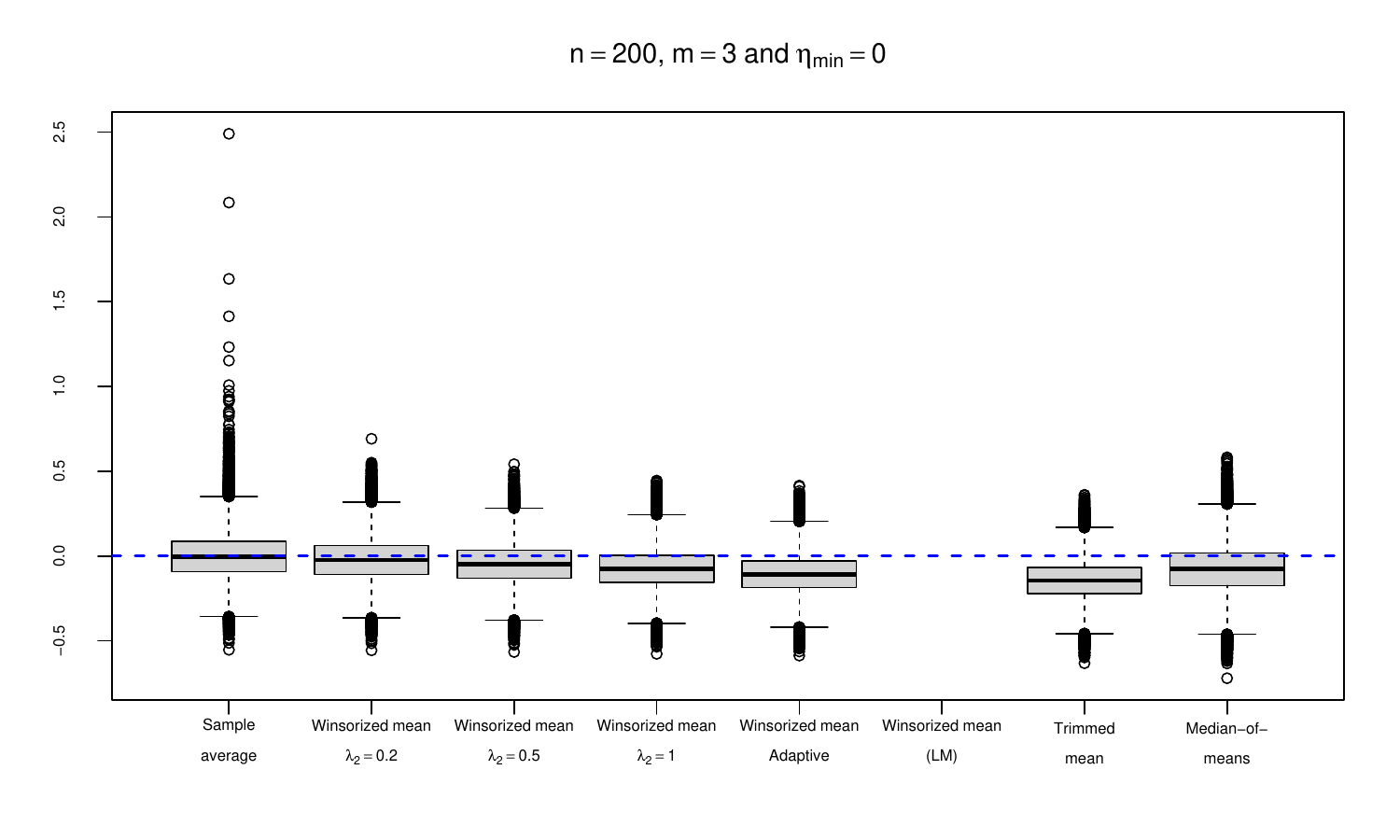}
		\includegraphics[width=10cm, height=7cm]{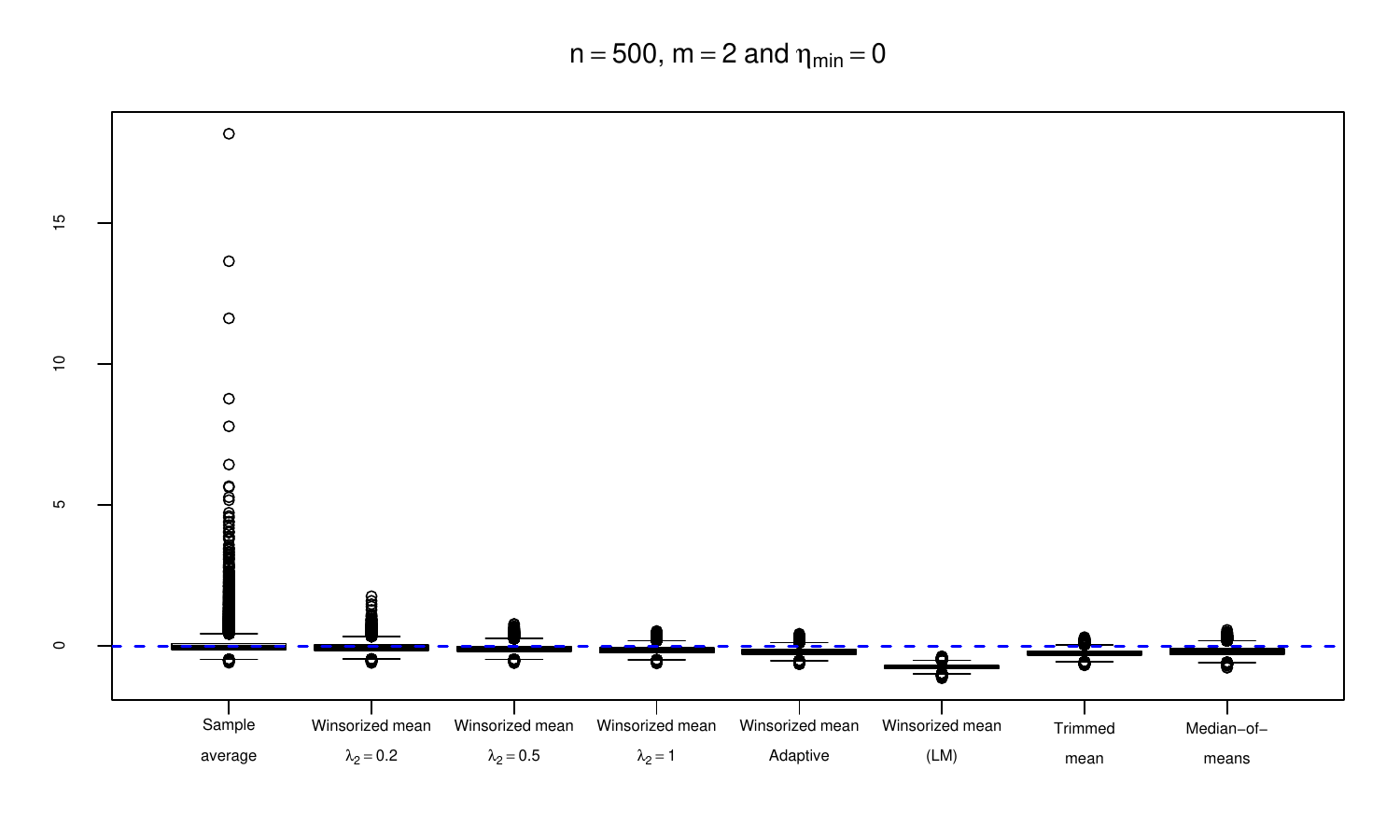}
	\includegraphics[width=10cm, height=7cm]{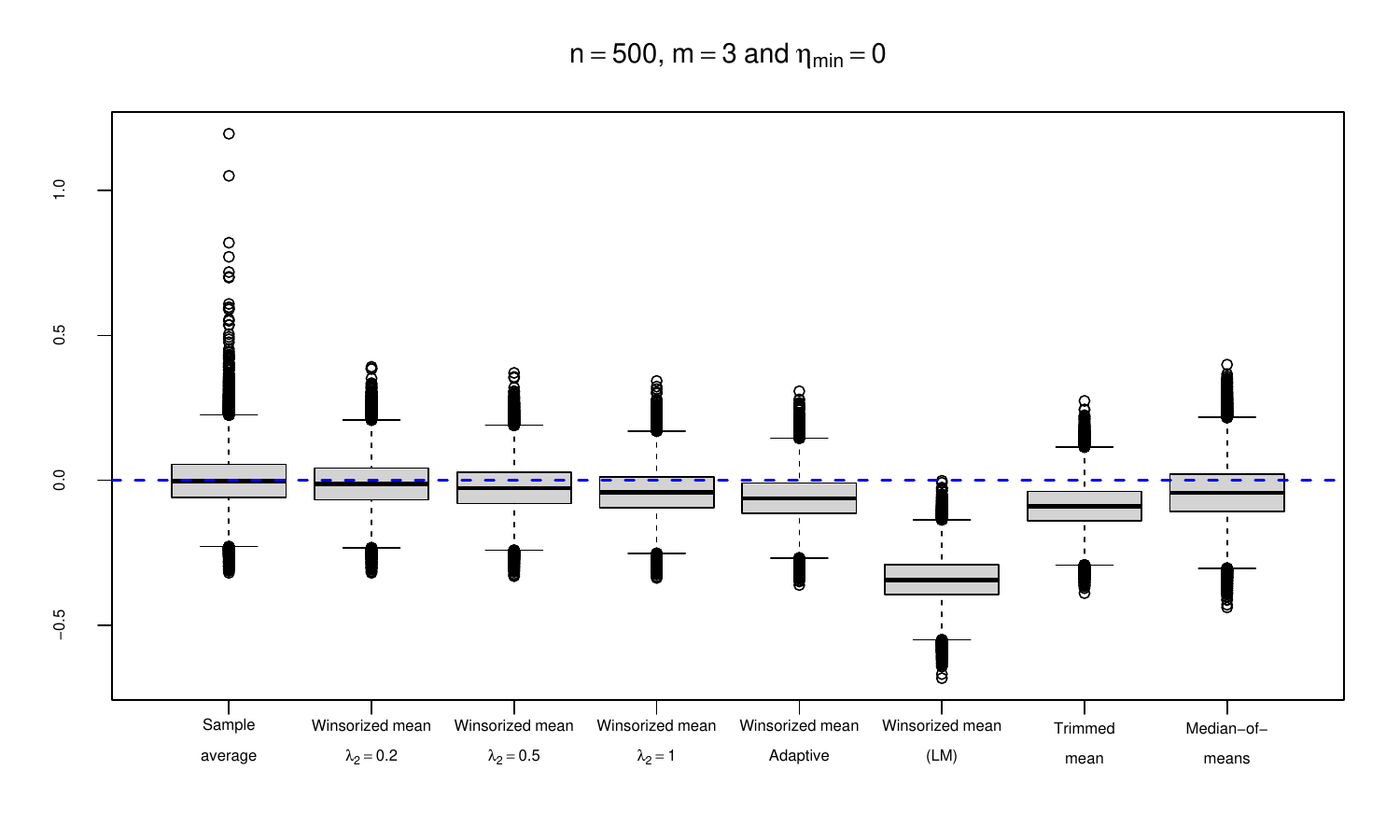}
	\caption{\footnotesize Box plots illustrating the distribution of the studied estimators. The dashed blue line indicates the true value (zero) of~$\mu$.}
	\label{fig:1}
\end{figure}
\end{landscape}

We first study a setting without contamination (i.e.,~$\eta_{\min}=0$). All non-adaptive estimators are implemented with~$\eta=0$. Table~\ref{tab:1} contains the mean absolute estimation errors whereas Figure~\ref{fig:1} contains box plots illustrating the distribution of the estimators.

As expected, the box plots reveal that the sample average has very heavy tails and can be rather erratic (in particular when~$m=2$). In implementing the winsorized mean,~$\lambda_2=0.2$ seems to work best, but the performance is not overly sensitive to the choice of~$\lambda_2$. 

In the numerical results, the adaptive winsorized mean estimator was implemented with~$\rho=0.9$ and~$\sigma_m$ being the oracle centered~$m$th absolute moment of the distribution of~$X_1$. Even with this choice, the adaptive winsorized mean estimator turned out to always pick~$\hat{g}=g_{\max}$. Furthermore, it turned out that~$\cap_{j=1}^{g_{\max}}\mathbb{I}(\eta_j)=\mathbb{B}\del[1]{\hat{\mu}_{n}(\eps_A(\eta_{g_{\max}})),B (\eta_{g_{\max}})}$, implying, by the definition of~$\hat{\mu}_{n,A}$, that~$\hat{\mu}_{n,A}=\hat{\mu}_{n}(\eps_A(\eta_{g_{\max}}))$. However, even though $\hat{\mu}_{n}(\eps_A(\eta_{g_{\max}}))$ uses the small~$0<\eta_{g_{\max}}=0.5\rho^{g_{\max}}\leq \log(6/\delta)/n$, it still winsorizes more observations than all of the~$\hat{\mu}_{n,\lambda_2}$ for~$\lambda_2\in\cbr[0]{0.2,0.5,1}$. This ``excessive'' winsorization explains its larger downward bias towards the median (which is negative).

Table~\ref{tab:1} shows that the mean absolute estimation error of the winsorized estimators is lower than that of the trimmed mean when~$X_i\sim \mathsf{m}_{t,\gamma}$. As mentioned, we also experimented with~$X_i\sim \mathsf{t}(\gamma)$ with~$\gamma\in\cbr[0]{2.01,3.01}$, for which the mean and median coincide. Here the winsorized and trimmed mean were both more precise than the sample average irrespective of the choice of~$\lambda_2$, but now the trimmed mean was slightly more precise than the winsorized mean. Since, e.g., Theorem 1.3.1 in~\cite{oliveira2025finite} establishes performance guarantees for the trimmed mean similar to those established for the winsorized mean in Theorem~\ref{thm:maintext}, it is not surprising that none of these two estimators uniformly dominates the other. 

The winsorized mean of~\cite{LM21} was not implementable for~$n=200$ as its~$\eps=24\log(4/\delta)/n>0.5$. When~$n=500$, their estimator is not very precise as it (essentially) uses~$\lambda_2=24$ and hence winsorizes so many observations that it approaches the median (which is negative). This underscores the importance for allowing ``small''~$\lambda_2$ as in our Theorem~\ref{thm:maintext}. 

\begin{table}[ht]
	
	\centering
	
	$\eta_{\min}=0$
	
	\medskip
	\begin{tabular}{rrrrrrrrr}
		\toprule
		& $S_n$ & $\hat{\mu}_{n,0.2}$ & $\hat{\mu}_{n,0.5}$ & $\hat{\mu}_{n,1}$ & $\hat{\mu}_{n,A}$ & $\hat{\mu}_{n,LM}$ & $\hat{\mu}_{n,T}$ & $\hat{\mu}_{n,MoM}$ \\ 
		\midrule
		\multirow{2}{*}{$n=200$\hspace{0.2cm}}
		$m=2$ & 0.224 & 0.199 & 0.215 & 0.257 & 0.314 & N/A & 0.379 & 0.318 \\ 
		$m=3$ & 0.106 & 0.103 & 0.106 & 0.114 & 0.130 & N/A & 0.157 & 0.133 \\ 
		\midrule
		\multirow{2}{*}{$n=500$\hspace{0.2cm}} 
		$m=2$ & 0.150 & 0.134 & 0.144 & 0.168 & 0.211 & 0.748 & 0.260 & 0.210\\
		$m=3$ & 0.068 & 0.066 & 0.067 & 0.071 & 0.080 & 0.343 & 0.098 & 0.085\\
		\bottomrule
	\end{tabular}
	\caption{\footnotesize Mean absolute estimation errors.~$S_n=n^{-1}\sum_{i=1}^nX_i$ denotes the sample average. $\hat{\mu}_{n,\lambda_2}=\hat{\mu}_n(\eps)=\hat{\mu}_n(\lambda_2\log(6/\delta)/n)$ denotes the winsorized mean estimator in~\eqref{eqn:winsmean} with~$\eps(\eta)$ chosen as in~\eqref{eq:epsfam}, which is always implemented with~$\lambda_1=1.01$ and with~$\lambda_2\in\cbr[0]{0.2,0.5,1}$. $\hat{\mu}_{n,A}$ is the adaptive estimator from Section~\ref{sec:adaptive}, which is always implemented with~$\lambda_1=1.5$ and $\lambda_2=0.2$. $\hat{\mu}_{n,LM}$ is the winsorized mean estimator from Section 2 in \cite{LM21}, $\hat{\mu}_{n,T}$ is the trimmed mean estimator from Theorem 1.3.1 in \cite{oliveira2025finite}, and $\hat{\mu}_{n,MoM}$ is the median-of-means estimator from Theorem 2 in \cite{lugosi2019mean}. }
	\label{tab:1}
\end{table}

\subsection{Contamination:~$\eta_{\min}=0.1$}
We next consider a setting where 10\% of the observations have been contaminated, amounting to~$\eta_{\min}=0.1$. All non-adaptive estimators are implemented with~$\eta=0.2$ to reflect that when there is contamination one  typically does not know the exact fraction of observations that have been contaminated. The adversary replaces~$0.1\cdot n$ randomly chosen observations by the 99th percentile of~$\mathsf{m}_{2,m+0.01}$. 

\begin{landscape}
\begin{figure}
	\includegraphics[width=10cm, height=7cm]{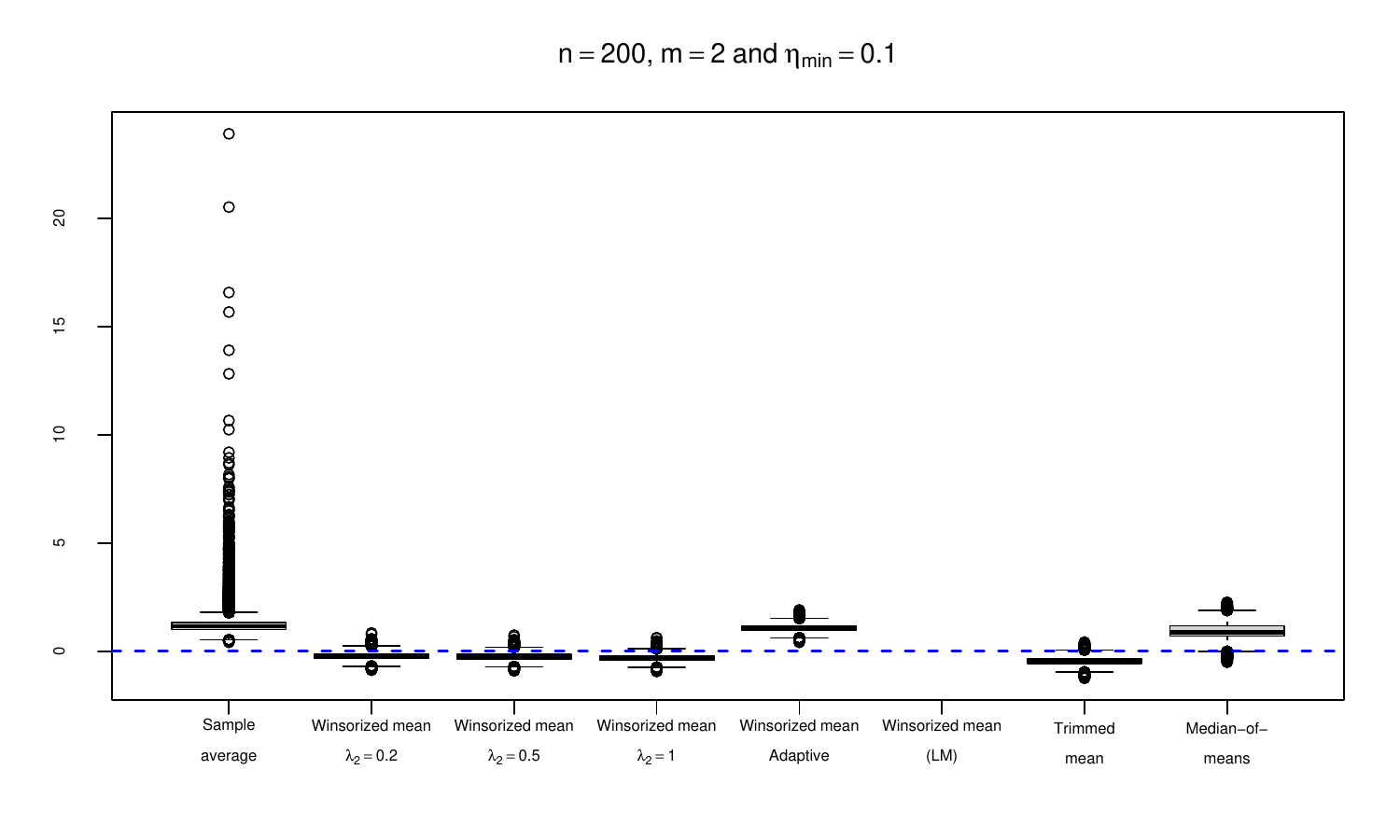}
	\includegraphics[width=10cm, height=7cm]{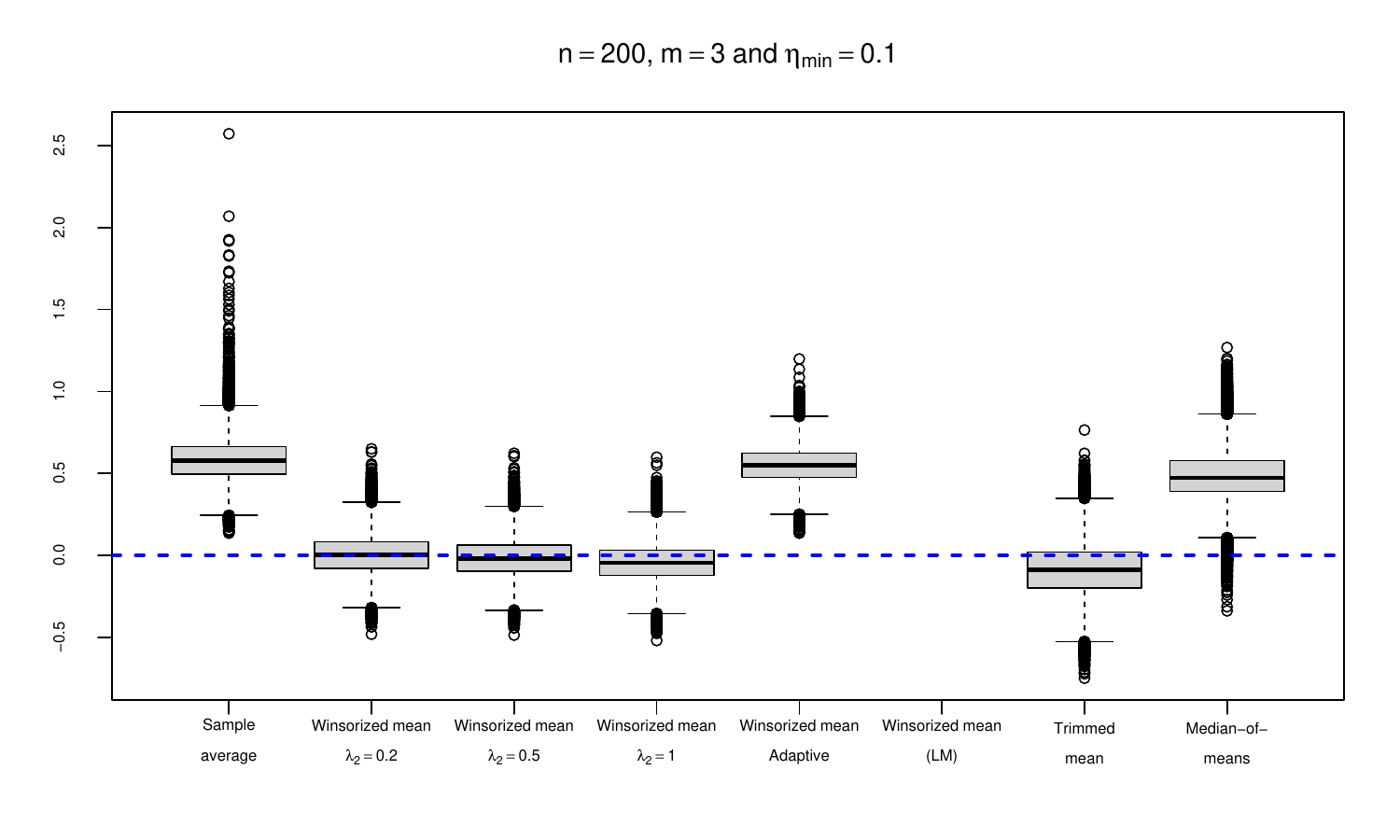}
		\includegraphics[width=10cm, height=7cm]{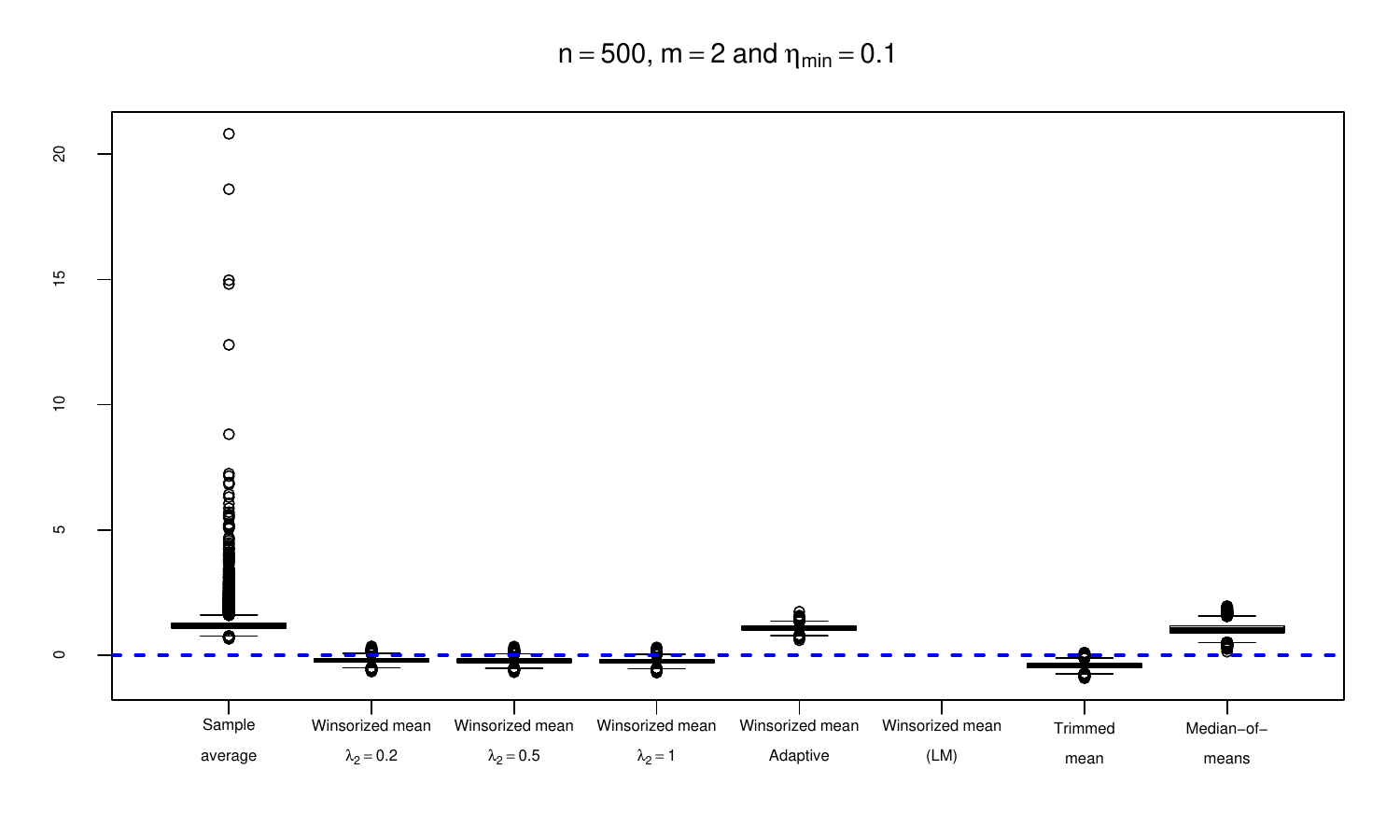}
	\includegraphics[width=10cm, height=7cm]{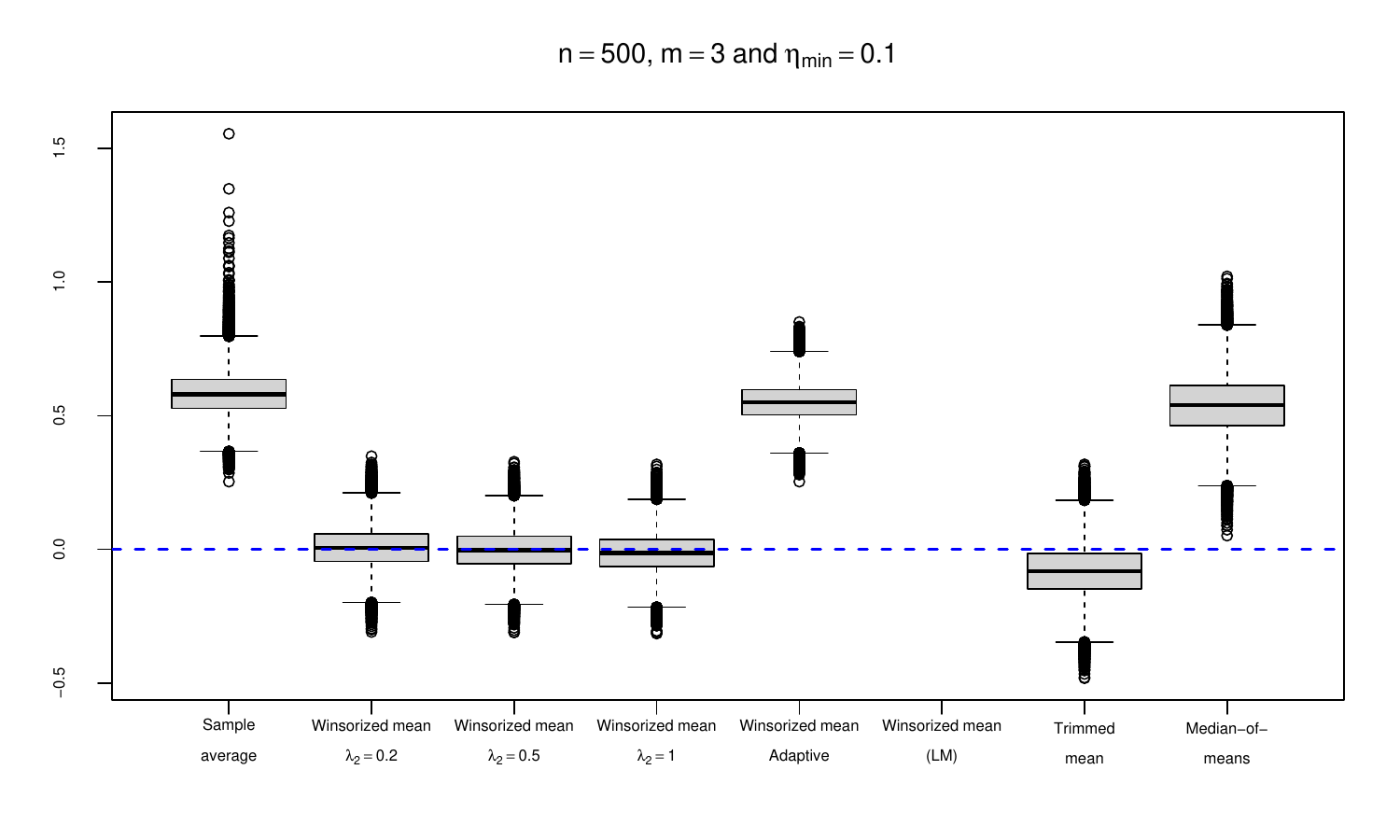}
	\caption{\footnotesize Box plots illustrating the distribution of the studied estimators. The dashed blue line indicates the true value (zero) of~$\mu$.}
	\label{fig:2}
\end{figure}
\end{landscape}

The mean absolute estimation errors can be found in Table~\ref{tab:2} and the box plots illustrating the distribution of the estimators can be found in Figure~\ref{fig:2}. The box plots reveal that despite contamination the distribution of the winsorized mean estimators from~\eqref{eqn:winsmean} with~$\eps(\eta)$ chosen as in~\eqref{eq:epsfam} is centered around the true mean irrespective of the value of~$m$ and~$n$. As explained already, the adaptive estimator~$\hat{\mu}_{n,A}$ frequently equals~$\hat{\mu}_{n}(\eps_A(\eta_{g_{\max}}))$. In the presence of contamination this means that ``too few'' observations are winsorized, explaining why it performs only slightly better than the sample average and is centered similarly.

The trimmed mean estimator has a larger downward bias than the winsorized mean estimators. However, when we implemented the winsorized and trimmed means with the ``oracle value''~$\eta=\eta_{\min}=0.1$ instead of~$\eta=0.2$, we found that the trimmed mean performed better than the winsorized mean (and the latter was most precise for~$\lambda_2=1$). As already discussed in the previous section, it is not surprising that neither of these estimators uniformly dominates the other. Finally, the winsorized mean estimator of~\cite{LM21} is not implementable as this requires~$\eta<1/16$.

\begin{table}[ht]
	
	\centering
	
	$\eta_{\min}=0.1$
	
	\medskip
	\begin{tabular}{rrrrrrrrr}
		\toprule
		& $S_n$ & $\hat{\mu}_{n,0.2}$ & $\hat{\mu}_{n,0.5}$ & $\hat{\mu}_{n,1}$ & $\hat{\mu}_{n,A}$ & $\hat{\mu}_{n,LM}$ & $\hat{\mu}_{n,T}$ & $\hat{\mu}_{n,MoM}$ \\ 
		\midrule
		
		\multirow{2}{*}{$n=200$\hspace{0.2cm}}
		$m=2$ & 1.202 & 0.237 & 0.266 & 0.311 & 1.076 & N/A & 0.446 & 0.902 \\ 
		$m=3$ & 0.583 & 0.096 & 0.095 & 0.100 & 0.550 & N/A & 0.149 & 0.482 \\ 
		\midrule
		\multirow{2}{*}{$n=500$\hspace{0.2cm}} 
		$m=2$ & 1.201 & 0.214 & 0.229 & 0.251 & 1.077 & N/A & 0.423 & 1.035\\
		$m=3$ & 0.583 & 0.061 & 0.060 & 0.061 & 0.551 & N/A & 0.104 & 0.540\\
		\bottomrule
	\end{tabular}
	\caption{\footnotesize Mean absolute estimation errors.~$S_n=n^{-1}\sum_{i=1}^n\tilde{X}_i$ denotes the sample average. $\hat{\mu}_{n,\lambda_2}=\hat{\mu}_n(\eps)=\hat{\mu}_n(1.01\cdot 0.2 +\lambda_2\log(6/\delta)/n)$ denotes the winsorized mean estimator in~\eqref{eqn:winsmean} with~$\eps(\eta)$ chosen as in~\eqref{eq:epsfam}, which is always implemented with~$\lambda_1=1.01$ and with~$\lambda_2\in\cbr[0]{0.2,0.5,1}$. $\hat{\mu}_{n,A}$ is the adaptive estimator from Section~\ref{sec:adaptive}, which is always implemented with~$\lambda_1=1.5$ and $\lambda_2=0.2$. $\hat{\mu}_{n,LM}$ is the winsorized mean estimator from Section 2 in \cite{LM21}, $\hat{\mu}_{n,T}$ is the trimmed mean estimator from Theorem 1.3.1 in \cite{oliveira2025finite}, and $\hat{\mu}_{n,MoM}$ is the median-of-means estimator from Theorem 2 in \cite{lugosi2019mean}. }
	\label{tab:2}
\end{table}

\bibliographystyle{ecta} 
\bibliography{ref}		

\newpage
\begin{appendix}
	\numberwithin{equation}{section}
	
	\section{Outline of the proof strategy for Theorem \ref{thm:maintext}}\label{sec:Decomposition}
	For~$p\in(0,1)$ and a random variable~$Z$, denote by~$Q_p(Z)$ the $p$-quantile of the distribution of~$Z$, that is
	\begin{equation}\label{eq:quantiledef}
		Q_p(Z)=\inf\cbr[1]{z\in \R:\P(Z\leq z)\geq p}.
	\end{equation}
	To prove Theorem \ref{thm:maintext}, we first establish in Lemma~\ref{lem:quantiles} (cf.~also Remark~\ref{rem:qlepschoice}) that on a set~$G_n$, say, of  probability at least~$1-\frac{4}{6}\delta$, one has that~$\hat{\alpha}=\tilde X_{\lceil \eps n \rceil}^*$ and~$\hat{\beta}=\tilde X_{\lfloor(1-\eps )n\rfloor+1}^*$ are bounded from above and below by suitable population quantiles:
	\begin{equation}\label{eq:lowerdef}
		Q_{c_1\eps}(X_1)=:\underline{\alpha}\leq \hat{\alpha}\leq \overline{\alpha}:= Q_{c_2\eps}(X_1),
	\end{equation}
	and
	\begin{equation}\label{eq:upperdef}
		Q_{1-c_2\eps}(X_1)=:\underline{\beta}\leq \hat{\beta}\leq \overline{\beta}:= Q_{1-c_1\eps}(X_1);
	\end{equation}
	here~$c_1\in(0,1)$,~$c_2\in(1,\infty)$ (cf.~Equations~\eqref{eq:c1} and~\eqref{eq:c2} for the precise definition of~$c_1$ and~$c_2$, respectively), and~$0 <\eps(c_1+c_2)<1$ holds, such that all expressions are well-defined. Together,~\eqref{eq:lowerdef} and~\eqref{eq:upperdef} imply, via obvious monotonicity properties of~$(a,b) \mapsto \phi_{a,b}$, that $$\phi_{\underline{\alpha},\underline{\beta}} \leq \phi_{\hat{\alpha},\hat\beta} \leq \phi_{\overline{\alpha},\overline{\beta}}.$$
	On~$G_n$ one thus obtains the following control of~$\frac{1}{n}\sum_{i=1}^n[\phi_{\hat{\alpha},\hat\beta}(\tilde{X}_i)-\mu]$:
	\begin{equation}\label{eq:lbub}
		\frac{1}{n}\sum_{i=1}^n\sbr[1]{\phi_{\underline{\alpha},\underline{\beta}}(\tilde{X}_i)-\mu}
		\leq
		\frac{1}{n}\sum_{i=1}^n\sbr[1]{\phi_{\hat{\alpha},\hat\beta}(\tilde{X}_i)-\mu}
		\leq
		\frac{1}{n}\sum_{i=1}^n\sbr[1]{\phi_{\overline{\alpha},\overline{\beta}}(\tilde{X}_i)-\mu}.
	\end{equation}
	Furthermore, the far right-hand side in~\eqref{eq:lbub} can be decomposed as 
	\begin{align}
		\frac{1}{n}\sum_{i=1}^n\sbr[1]{\phi_{\overline{\alpha},\overline{\beta}}(\tilde{X}_i)-\mu}
		&=
		\underbrace{\frac{1}{n}\sum_{i=1}^n\sbr[1]{\phi_{\overline{\alpha},\overline{\beta}}(\tilde{X}_i)-\phi_{\overline{\alpha},\overline{\beta}}(X_i)}}_{\overline I_{n,1}} \notag\\
		&+
		\underbrace{\frac{1}{n}\sum_{i=1}^n\sbr[1]{\phi_{\overline{\alpha},\overline{\beta}}(X_i)-\E\phi_{\overline{\alpha},\overline{\beta}}(X_i)}}_{\overline I_{n,2}}\notag\\
		&+\underbrace{\frac{1}{n}\sum_{i=1}^n\sbr[1]{\E\phi_{\overline{\alpha},\overline{\beta}}(X_i)-\mu}}_{\overline I_{n,3}}\label{eq:decomp}.
	\end{align}
	Thus, it suffices to control:
	\begin{enumerate}
		\item $\overline I_{n,1}$, i.e.,~an error incurred from computing the winsorized mean on the corrupted data $\tilde{X}_1,\hdots,\tilde{X}_n$ instead of the uncorrupted~$X_1,\hdots,X_n$;
		\item $\overline I_{n,2}$, i.e.,~the difference between the sample and population means of the bounded~$\phi_{\overline\alpha,\overline\beta}$ evaluated at the uncorrupted data; and
		\item $\overline I_{n,3}$, i.e., a difference between the winsorized and raw population means.
	\end{enumerate}
	Replacing $\phi_{\overline{\alpha},\overline{\beta}}$ by~$\phi_{\underline{\alpha},\underline{\beta}}$ in~$\overline{I}_{n,k}$ for~$k = 1, 2, 3$ and denoting the obtained quantities~$\underline{I}_{n,k}$ for~$k = 1, 2, 3$, the left-hand side of~\eqref{eq:lbub} can be decomposed analogously as
	\begin{equation}\label{eq:decomp2}
		\frac{1}{n}\sum_{i=1}^n\sbr[1]{\phi_{\underline{\alpha},\underline{\beta}}(\tilde{X}_i)-\mu}
		=
		\underline I_{n,1}+\underline I_{n,2}+\underline I_{n,3}.
	\end{equation}
	Lemmas \ref{lem:noisecontrol},~\ref{lem:trimmeanconc}, and \ref{lem:meancontrol} in Section~\ref{sec:Is} are auxiliary results that allow us to bound the~$\underline{I}_{n,i}$ and~$\overline{I}_{n,i}$. The proof of Theorem~\ref{thm:maintext} collects the respective expressions and concludes. 
	
	\section{Some preparatory lemmas}\label{sec:lemmas}
	The functions~$h_+:[0,\infty)\to [0,\infty)$ and~$h_{-}:[0,1)\to [0,\infty)$ defined as
	\begin{equation}\label{eq:hs}
		h_+(\nu):=(1+\nu)\log(1+\nu)-\nu\qquad\text{and}\qquad h_-(\nu):=(1-\nu)\log(1-\nu)+\nu
	\end{equation}
	will enter in the following lemmas. 
	
	We first recall suitable versions of the classic lower and upper multiplicative Chernoff bounds for the Bernoulli distribution from~\cite{hagerup1990guided}. The first is taken from their Equation~(5), and the second from the equation preceding their Equation~(7).
	\begin{lemma}\label{lem:Chernoff}
		Let~$B$ be binomially distributed with success probability~$p\in(0,1)$ and number of trials~$n\in\N$. Then 
		\begin{enumerate}
			\item $\P\del[1]{B\geq (1+\nu)np}
			\leq
			e^{-nph_+(\nu)}$ for every $\nu\in(0,\infty)$.
			\item $\P\del[1]{B\leq  (1-\nu)np}
			\leq
			e^{-nph_-(\nu)}$ for every $\nu\in(0,1)$.
		\end{enumerate}
	\end{lemma}
	
	The following lemma and its proof make use of some elementary properties of Lambert's~$W$ function (cf., e.g., \cite{Corless1996}).
	
	\begin{lemma}\label{lem:Prophfs}
		For given~$\lambda_1\in(1,\infty)$ and~$\eta \in [0, 1]$, we make the following observations.
		\begin{enumerate}
			\item Define~$A_+:=1-\lambda_1^{-1}\mathds{1}\cbr[0]{\eta>0}$,~$\nu_{+}(c):=\frac{A_+}{c}-1$, and~$f(c):=ch_+(\nu_{+}(c))$ for~$c\in(0,A_+)$. Then,
			\begin{enumerate}
				\item $f$ is differentiable and strictly decreasing on~$(0,A_+)$, and
				\item $\lim_{c\downarrow 0}f(c)=\infty$ and~$\lim_{c\uparrow A_+}f(c)=0$.
			\end{enumerate}
			In particular,~$f$ is a bijection from~$(0,A_+)$ to~$(0,\infty)$ with inverse
			\begin{equation}\label{eq:eq+}
				f^{-1}(r)=-A_+W_0(-e^{-(r+A_+)/A_+}),
			\end{equation}
			where~$W_0$ is the principal branch of Lambert's~$W$ function, and
			\begin{equation}\label{eq:+bounds}
				A_+e^{-(r+A_+)/A_+}\leq f^{-1}(r) < A_+.
			\end{equation}

			\item Define~$A_-:=1+\lambda_1^{-1}\mathds{1}\cbr[0]{\eta>0}$,~$\nu_-(c):=1-\frac{A_-}{c}$, and~$g(c):=ch_-(\nu_-(c))$ for~$c\in(A_-,\infty)$. Then,
			\begin{enumerate}
				\item $g$ is differentiable and strictly increasing on~$(A_-,\infty)$, and
				\item $\lim_{c\downarrow A_-}g(c)=0$ and~$\lim_{c\uparrow \infty}g(c)=\infty$.
			\end{enumerate}
			In particular,~$g$ is a bijection from~$(A_-,\infty)$ to~$(0,\infty)$ with inverse
			\begin{equation}\label{eq:eq-}
				g^{-1}(r)=-A_-W_{-1}(-e^{-(r+A_-)/A_-}),
			\end{equation}	
			where~$W_{-1}$ is the lower branch of Lambert's~$W$ function, and
			\begin{equation}\label{eq:-bounds}
				A_-+r\leq g^{-1}(r) \leq A_-+r+\sqrt{r^2+2A_-r}.
			\end{equation}
		\end{enumerate}
	\end{lemma}
	
	\begin{proof}
		Concerning Part 1., because the image of~$(0,A_+)$ under $\nu_+$ is~$(0,\infty)$, which is a subset of the domain of~$h_+$, it follows that~$f$ is well-defined. Next, note that
		\begin{equation}\label{eq:f}
			f(c)
			=
			ch_+(\nu_{+}(c))
			=
			A_+\log\del[2]{\frac{A_+}{c}}+c-A_+.
		\end{equation}
		Thus,~$f'(c)=1-A_+/c<0$ for~$c\in(0,A_+)$,~such that~$f$ is strictly decreasing. It also follows that $\lim_{c\downarrow 0}f(c)=\infty$ and $\lim_{c\uparrow A_+}f(c)=0$. As a consequence,~$f: (0, A_+) \to (0, \infty)$ has an inverse~$f^{-1}$, say. Fix an arbitrary~$r \in (0, \infty)$. Abbreviating~$z_r:=f^{-1}(r)/A_+$ and~$C_r:=r/A_+$, it follows from~\eqref{eq:f} applied to~$c = f^{-1}(r)$ that
		\begin{equation}\label{eq:faux}
			z_r-1-\log(z_r)=C_r	\qquad\Longleftrightarrow \qquad e^{-z_r}(-z_r)=-e^{-(C_r+1)}.
		\end{equation}
		Noting that~$-e^{-(C_r+1)}\in(-e^{-1},0)$, we conclude that\footnote{Since~$-e^{-(C_r+1)}\in (-e^{-1},0)$, there are two real~$u$ solving~$e^uu=-e^{-(C_r+1)}$, which can be expressed in terms of the principal and lower branch of Lambert's~$W$ function, respectively. However, only the principal branch results in~$f^{-1}(r) \in(0,A_+)$.} 
		\begin{equation*}
			-z_r=W_0(-e^{-(C_r+1)}) \qquad \Leftrightarrow \qquad f^{-1}(r)=-A_+W_0(-e^{-(r+A_+)/A_+})\in(0,A_+).
		\end{equation*}
		The claimed lower bound on~$f^{-1}(r)$ follows from~\eqref{eq:faux}, since~$z_r\in(0,1)$ such that
		\begin{equation*}
			e^{-z_r}(-z_r)=-e^{-(C_r+1)} \quad \Rightarrow \quad z_r\geq e^{-(C_r+1)} \quad \Leftrightarrow \quad f^{-1}(r)\geq A_+e^{-(r+A_+)/A_+}.
		\end{equation*}

		Concerning Part~2., because the image of~$(A_-,\infty)$ under~$\nu_-$ is~$(0,1)$, which is a subset of the domain of~$h_-$, it follows that~$g$ is well-defined. Next, note that
		\begin{equation}\label{eq:g}
			g(c)
			=
			ch_-(\nu_{-}(c))
			=
			A_-\log\del[2]{\frac{A_-}{c}}+c-A_-.
		\end{equation}
		Thus,~$g'(c)=1-A_-/c>0$ for~$c\in(A_-,\infty)$,~such that~$g$ is strictly increasing. It also follows that $\lim_{c\downarrow A_-}g(c)=0$ and
		\begin{equation*}
			\lim_{c\uparrow \infty}g(c)
			=
			\lim_{c\uparrow \infty}c\cdot\del[2]{\frac{A_-\log(A_-)}{c}-\frac{A_-\log(c)}{c}+1-\frac{A_-}{c}}=\infty.
		\end{equation*}
		As a consequence,~$g: (A_-, \infty) \to (0, \infty)$ has an inverse~$g^{-1}$, say. Fix an arbitrary~$r \in (0, \infty)$. Re-defining~$z_r:=g^{-1}(r)/A_-$ and~$C_r:=r/A_-$, it follows from~\eqref{eq:g} applied to~$c = g^{-1}(r)$ that
		\begin{equation}\label{eq:gaux}
			z_r-1-\log(z_r)=C_r	\qquad\Longleftrightarrow \qquad e^{-z_r}(-z_r)=-e^{-(C_r+1)}.
		\end{equation}
		With the new definitions of~$z_r$ and~$C_r$ in place, the display~\eqref{eq:gaux} is identical to~\eqref{eq:faux}. Thus, arguing as after~\eqref{eq:faux}, it follows that
		\begin{equation*}
			g^{-1}(r)=-A_-W_{-1}(-e^{-(r+A_-)/A_-})\in(A_-,\infty);	
		\end{equation*} 
		where we note that it is now only the \emph{lower} branch of Lambert's~$W$ function that results in~$g^{-1}(r)\in(A_-,\infty)$. The claimed lower bound on~$g^{-1}(r)$ follows from~\eqref{eq:gaux} since~$z_r\in(1,\infty)$ such that
		\begin{equation*}
			z_r-1\geq z_r-1-\log(z_r)=C_r	\quad \Longleftrightarrow \quad z_r\geq C_r+1
			\quad \Longleftrightarrow \quad g^{-1}(r) \geq r+A_-.
		\end{equation*}
		Next, to provide the claimed upper bound on~$g^{-1}(r)$, recall the standard inequality~$$\log(z)\leq z-1-(z-1)^2/(2z) \quad \text{ for } z\geq 1,$$ which used in~\eqref{eq:gaux} implies that 
		\begin{equation*}
			\frac{(z_r-1)^2}{2z_r}
			\leq
			z_r-1-\log(z_r)
			=
			C_r
			\qquad
			\Longrightarrow
			\qquad 
			z_r^2-2(1+C_r)z_r+1\leq 0.
		\end{equation*}
		Noting that the coefficient on~$z_r^2$ is positive, solving for the roots of this second degree polynomial yields that~$z_r\leq 1+C_r+\sqrt{C_r(C_r+2)}$. Therefore, recalling that~$z_r=g^{-1}(r)/A_-$ and~$C_r=r/A_-$, one concludes that~$g^{-1}(r)
		\leq 
		A_-+r+\sqrt{r^2+2A_-r}.$
	\end{proof}
	Recall the notation of Lemma~\ref{lem:Prophfs} (in particular~$f^{-1}$ and~$g^{-1}$,~$A_+=1-\lambda_1^{-1}\mathds{1}\cbr[0]{\eta>0}$, and~$A_-=1+\lambda_1^{-1}\mathds{1}\cbr[0]{\eta>0}$), and \emph{throughout the remainder of the paper} define, for every~$\epsilon \in (0, \infty)$ and~$\delta \in (0, \infty)$, the quantities
	\begin{equation}\label{eq:c1}
		c_1 := f^{-1}\left(\log(6/\delta)/(n\epsilon)\right) = -A_+W_0\del[1]{-e^{-(\frac{\log(6/\delta)}{\epsilon n}+A_+)/A_+}}\in (0,A_+),
	\end{equation}
	as well as
	\begin{equation}\label{eq:c2}
		c_2 := g^{-1}\left(\log(6/\delta)/(n\epsilon)\right) = -A_-W_{-1}\del[1]{-e^{-(\frac{\log(6/\delta)}{\epsilon n}+A_-)/A_-}}\in(A_-,\infty).
	\end{equation}
	We emphasize that in addition to~$\epsilon, n$ and~$\delta$, the quantities~$c_1$ and~$c_2$ also depend on~$\lambda_1$ and~$\eta$, although none of these dependencies is shown explicitly. Despite these dependencies, the following lemma (which is written with applications to the case~$\epsilon = \varepsilon$ as in~\eqref{eq:epsfam} in  mind, but applies more generally) bounds~$c_1$ and~$c_2$ in terms of the parameters~$\lambda_1$ and~$\lambda_2$ only. 
	\begin{lemma}\label{lem:cControl}
		Let~$n \in \N$,~$\delta \in (0, 1)$,~$\lambda_1 \in (1, \infty)$,~$\lambda_2 \in (0, \infty)$,~$\eta \in [0, 1]$, and suppose that~$\epsilon \in (0, 1)$ satisfies~$\epsilon \geq \lambda_2 \log(6/\delta)/n$. Then, for~$c_1$ as defined in~\eqref{eq:c1} and~$c_2$ as defined in~\eqref{eq:c2}, it holds that
		\begin{align*}
			0<(1-\lambda_1^{-1})\exp \del[2]{{-\frac{1}{\lambda_2(1-\lambda_1^{-1})}-1}}\leq\ &c_1 < A_+ \leq 1,\\
			1 \leq A_- <\ &c_2 \leq 2+\lambda_2^{-1}+\sqrt{\lambda_2^{-2}+4\lambda_2^{-1}},
		\end{align*}
		and that
		\begin{align}\label{eq:welldefined}
			0
			&<
			\epsilon\min(c_1,c_2) \notag\\
			&\leq 
			\epsilon(c_1+c_2)\notag\\
			&\leq 
			2\epsilon +\frac{\log(6/\delta)}{n}+\sqrt{\del[2]{\frac{\log(6/\delta)}{n}}^2+2\sbr[1]{1+\lambda_1^{-1}\mathds{1}(\eta>0)}\frac{\log(6/\delta)}{n}\epsilon}\notag\\
			&\leq
			2\epsilon +\frac{\log(6/\delta)}{n}+\sqrt{\del[2]{\frac{\log(6/\delta)}{n}}^2+4\frac{\log(6/\delta)}{n}\epsilon}.
		\end{align}	
	\end{lemma}
	\begin{proof}
		Throughout this proof set~$r :=\log(6/\delta)/(n\epsilon)\leq \lambda_2^{-1}$. Note that~$c_1 = f^{-1}\left(r\right) <A_+\leq 1$. The lower bound in~\eqref{eq:+bounds}, using~$A_+=1-\lambda_1^{-1}\mathds{1}(\eta>0)\geq 1-\lambda_1^{-1}>0$ since~$\lambda_1^{-1}<1$, yields
		\begin{equation*}
			c_1\geq A_+e^{-(r+A_+)/A_+}
			=
			A_+e^{-r/A_+-1}
			\geq
			(1-\lambda_1^{-1})\exp \del[2]{{-\frac{1}{\lambda_2(1-\lambda_1^{-1})}-1}}
			>0.
		\end{equation*}	
		Similarly, since~$c_2 = g^{-1}(r) > A_- \geq 1$, the upper bound in~\eqref{eq:-bounds}, using~$A_-=1+\lambda_1^{-1}\mathds{1}(\eta>0)\leq 2$, yields
		\begin{equation*}
			c_2\leq A_-+r+\sqrt{r^2+2A_-r}
			\leq 
			2+\lambda_2^{-1}+\sqrt{\lambda_2^{-2}+4\lambda_2^{-1}}.
		\end{equation*}

		Finally, since~$c_1$ and~$c_2$ are both strictly positive, it follows that~$$0 < \epsilon\min(c_1,c_2)\leq \epsilon(c_1+c_2),$$ and by~\eqref{eq:+bounds} and \eqref{eq:-bounds} of Lemma~\ref{lem:Prophfs}, as well as similar reasoning as above, 
		\begin{align*}
			\epsilon(c_1+c_2)
			&\leq
			\epsilon\del[1]{A_++A_-+r+\sqrt{r^2+2A_-r}}\\
			&= 
			\epsilon\del[4]{2+\frac{\log(6/\delta)}{n\epsilon}+\sqrt{\del[2]{\frac{\log(6/\delta)}{n\epsilon}}^2+2\sbr[1]{1+\lambda_1^{-1}\mathds{1}(\eta>0)}\frac{\log(6/\delta)}{n\epsilon}}}\\
			&=2\epsilon +\frac{\log(6/\delta)}{n}+\sqrt{\del[2]{\frac{\log(6/\delta)}{n}}^2+2\sbr[1]{1+\lambda_1^{-1}\mathds{1}(\eta>0)}\frac{\log(6/\delta)}{n}\epsilon},
		\end{align*}
		from which~\eqref{eq:welldefined} follows because~$\lambda_1^{-1}<1$.
	\end{proof}

	The following auxiliary lemma allows us to impose in the proof of Lemma~\ref{lem:quantiles} below (without loss of generality) the additional condition that the cdf of the~$X_i$ is continuous.
	\begin{lemma}\label{lem:contdisc}
		Fix~$n \in \N$ and~$\eta \in [0, 1]$. Suppose the numbers~$a \in \N \cap [1, n]$,~$b \in (0, 1)$, and~$\rho \in [0, 1]$ are such that\footnote{We denote by~$(\Omega, \mathcal{A}, \mathbb{P})$ the probability space on which the random variables $X_1, \hdots, X_n$ and $\tilde{X}_1, \hdots, \tilde{X}_n$ are defined.}
		\begin{equation}\label{eqn:indelcont}
			\P\left(\tilde{X}^*_a \geq Q_b(X_1)\right)\geq \rho,
		\end{equation}
		whenever the following conditions are satisfied:
		\begin{enumerate}[label=(\roman*)]
			\item $X_1, \hdots, X_n$ are i.i.d.~random variables,
			\item the random variables~$X_1, \hdots, X_n$ and~$\tilde{X}_1, \hdots, \tilde{X}_n$ satisfy~\eqref{eq:contamfrac}, and
			\item the cdf of~$X_1$ is continuous.
		\end{enumerate}   
		Then, whenever~(i)~and (ii) (but not necessarily (iii)) are satisfied, we have
		\begin{equation}\label{eqn:indelcontstar}
			\P\left(\tilde{X}^*_a \geq Q_b(X_1)\right) \geq \rho \quad \text{ and } \quad \P\left(-\tilde{X}^*_{n-a+1} \geq Q_{b}(-X_1)\right) \geq \rho.
		\end{equation}
		If all three inequality signs inside the probabilities in~\eqref{eqn:indelcont} and~\eqref{eqn:indelcontstar} are changed from ``$\geq$'' to~``$\leq$'', respectively, then the so-obtained statement is correct. 
	\end{lemma}
	
	\begin{proof}
		Fix~$n$ and~$\eta$ as in the first sentence of Lemma~\ref{lem:contdisc}, and suppose that (for the given numbers~$a, b$ and~$\rho$) the second sentence in Lemma~\ref{lem:contdisc} is a correct statement. Suppose that~$X_1, \hdots, X_n$ and~$\tilde{X}_1, \hdots, \tilde{X}_n$ satisfy (i)~and (ii) in Lemma~\ref{lem:contdisc} (but not necessarily satisfy (iii)). We show that then~\eqref{eqn:indelcontstar} holds. To this end, let~$U_i$ for~$i = 1, \hdots, n$ be independent, uniformly  distributed random variables on~$[-1, 1]$, that are independent of~$X_1, \hdots, X_n$ and~$\tilde{X}_1, \hdots, \tilde{X}_n$.\footnote{Such random variables~$U_1, \hdots, U_n$ certainly exist after suitably enlarging the probability space~$(\Omega, \mathcal{A}, \mathbb{P})$ on which~$X_1, \hdots, X_n$ and~$\tilde{X}_1, \hdots, \tilde{X}_n$ are defined. We don't spell out this (standard) enlargement argument for simplicity of notation, and assume without loss of generality that the~$U_i$ as required already exist on~$(\Omega, \mathcal{A}, \mathbb{P})$.} Fix~$k \in \N$, and define~$Y_{i,k}:=X_i+U_i/k$ for~$i = 1, \hdots, n$, which are i.i.d.~random variables. Because~$U_1$ has a continuous cdf, also~$Y_{1,k}$ has a continuous cdf (which can be shown by, e.g., combining Tonelli's theorem and the Dominated Convergence Theorem). Setting~$\tilde{Y}_{i,k}:=\tilde{X}_i+U_i/k$ for~$i = 1, \hdots, n$, we note that~$Y_{i,k} = \tilde{Y}_{i,k}$ is equivalent to~$X_i = \tilde{X}_i$, so that the random variables~$Y_{1,k}, \hdots, Y_{n,k}$ and~$\tilde{Y}_{1,k}, \hdots, \tilde{Y}_{n,k}$ satisfy~\eqref{eq:contamfrac}. The statement formulated in the second sentence of Lemma~\ref{lem:contdisc} is therefore applicable to~$Y_{1,k}, \hdots, Y_{n,k}$ and~$\tilde{Y}_{1,k}, \hdots, \tilde{Y}_{n,k}$, and delivers
		\begin{equation}\label{eqn:conclcont}
			\P\left(\tilde{Y}^*_{a,k} \geq Q_b(Y_{1,k})\right) \geq \rho.
		\end{equation}
		From~$X_1 - k^{-1} \leq Y_{1,k} \leq X_1 + k^{-1}$ and elementary equivariance and monotonicity properties of the map~$Q_p(\cdot)$ (defined in~\eqref{eq:quantiledef}), it follows that 
		\begin{equation}\label{qintcont}
			Q_{p}(X_1)-k^{-1} \leq Q_p(Y_{1,k}) \leq Q_{p}(X_1)+k^{-1} \quad \text{ for every } p \in (0, 1).
		\end{equation}
		From~$\tilde{Y}_{i,k}\leq \tilde{X}_i+k^{-1}$ for~$i = 1, \hdots, n$,  we obtain~$\tilde Y_{a,k}^*\leq \tilde X_{a}^*+k^{-1}$. Thus, whenever~$\tilde{Y}^*_{a,k} \geq Q_b(Y_{1,k})$, we have 
		\begin{equation*}
			\tilde X_{a}^*
			\geq
			\tilde Y_{a,k}^*-k^{-1}
			\geq
			Q_{b}(Y_{1,k})-k^{-1}
			\geq
			Q_{b}(X_1)-2k^{-1}.
		\end{equation*}
		Together with~\eqref{eqn:conclcont} we can conclude that~$	\P(\tilde X_{a}^* \geq Q_{b}(X_1)-2k^{-1}) \geq \rho$. Because~$k \in \N$ was arbitrary, we hence obtain the first inequality in~\eqref{eqn:indelcontstar} from 
		\begin{align*}
			\P(\tilde X_{a}^*\geq Q_{b}(X_1))
			&=
			\P\left(\bigcap_{k=1}^\infty \{\tilde X_{a}^* \geq Q_{b}(X_1)-2k^{-1}\}\right) \\
			&=
			\lim_{k\to\infty}\P(\tilde X_{a}^* \geq Q_{b}(X_1)-2k^{-1})
			\geq
			\rho.
		\end{align*}
		Summarizing, we have shown that~$\P(\tilde{X}^*_a \geq Q_b(X_1))\geq \rho$ whenever~$X_1, \hdots, X_n$ and~$\tilde{X}_1, \hdots, \tilde{X}_n$ satisfy (i)~and (ii). Note that~$X_1, \hdots, X_n$ and~$\tilde{X}_1, \hdots, \tilde{X}_n$ satisfy (i)~and (ii), if and only if~$-X_1, \hdots, -X_n$ and~$-\tilde{X}_1, \hdots, -\tilde{X}_n$ satisfy (i)~and (ii). We can hence apply the already established statement also to~$-X_1, \hdots, -X_n$ and~$-\tilde{X}_1, \hdots, -\tilde{X}_n$ to conclude~$\P((-\tilde{X})^*_a \geq Q_b(-X_1))\geq \rho$. Because~$-\tilde{X}^*_{n-a+1} = (-\tilde{X})^*_a$, the statement~$\P((-\tilde{X})^*_a \geq Q_b(-X_1))\geq \rho$ is equivalent to~$\P(-\tilde{X}^*_{n-a+1} \geq Q_{b}(-X_1)) \geq \rho$, so that we are done.
		
		To prove the remaining statement, we can use the same argument and construction as  that leading up to~\eqref{eqn:conclcont}, but now conclude~$\P(\tilde{Y}^*_{a,k} \leq Q_b(Y_{1,k})) \geq \rho$. From~$\tilde{Y}_{i,k} \geq \tilde{X}_i - k^{-1}$ for~$i = 1, \hdots, n$, we obtain~$\tilde{Y}_{a,k}^* \geq \tilde{X}_a^* - k^{-1}$. Thus, whenever~$\tilde{Y}^*_{a,k} \leq Q_b(Y_{1,k})$, we have (recall~\eqref{qintcont})
		\begin{equation}
			\tilde{X}^*_a \leq \tilde{Y}_{a,k}^* + k^{-1} \leq Q_b(Y_{1,k}) + k^{-1} \leq Q_b(X_1) + 2k^{-1}.
		\end{equation}
		Hence, under the condition that~$\P(\tilde{Y}^*_{a,k} \leq Q_b(Y_{1,k})) \geq \rho$, we obtain~$\P(\tilde{X}_a^* \leq Q_b (X_1) + 2k^{-1}) \geq \rho$. Because~$k \in \N$ was arbitrary, we can therefore conclude that
		\begin{equation}
			\P(\tilde{X}_a^* \leq Q_b (X_1) ) = \lim_{k \to \infty} \P\left(\tilde{X}_a^* \leq Q_b (X_1) + 2k^{-1}\right)  \geq \rho.
		\end{equation}
		Arguing as in the previous paragraph  establishes~$\P(-\tilde{X}^*_{n-a+1} \leq Q_{b}(-X_1)) \geq \rho$.
	\end{proof}

	The following lemma shows that (certain) order statistics of the contaminated data are close to related population quantiles of the uncontaminated data.
	\begin{lemma}\label{lem:quantiles}
		Let~$n\in\N$,~$\delta\in(0,1)$,~$\lambda_1 \in (1, \infty)$,~and~$\eta \in [0, 1]$. Let~$X_1,\hdots,X_n$ be i.i.d., and~\eqref{eq:contamfrac} be satisfied. Recall~$c_1$ from~\eqref{eq:c1} and~$c_2$ from~\eqref{eq:c2}, and let~$\epsilon \in (0, 1)$ satisfy
		\begin{equation}\label{eqn:condeps}
			\epsilon \geq \lambda_1 \eta \quad \text{ and } \quad \epsilon c_2 <1.
		\end{equation}
		Then, each of \eqref{eq:LBLQ}--\eqref{eq:UBUQ} below holds with probability at least~$1-\delta/6$:
		\begin{eqnarray}
			\tilde X_{\lceil \epsilon n \rceil}^* 
			&\geq&
			Q_{c_1\epsilon}(X_1)\label{eq:LBLQ};\\
			\tilde X_{\lceil(1-\epsilon)n\rceil}^*
			&\geq&
			Q_{1-c_2\epsilon}(X_1)\label{eq:LBUQ};\\
			\tilde{X}_{\lfloor \epsilon n \rfloor+1}^*
			&\leq&
			Q_{c_2\epsilon}(X_1)\label{eq:UBLQ};\\
			\tilde{X}_{\lfloor (1-\epsilon)n\rfloor+1}^*
			&\leq&
			Q_{1-c_1\epsilon}(X_1).\label{eq:UBUQ}
		\end{eqnarray}	
	\end{lemma}	
	
	\begin{remark}
		Inspection of the proof of Lemma~\ref{lem:quantiles} shows that one does \emph{not} need to impose the condition~$\epsilon c_2<1$ in~\eqref{eqn:condeps} to establish only the probability statements concerning the inequalities in~\eqref{eq:LBLQ} and~\eqref{eq:UBUQ}. 
	\end{remark}
	
	\begin{remark}\label{rem:qlepschoice}
		The conditions~$\epsilon \in (0, 1)$ and~\eqref{eqn:condeps} are satisfied for~$\epsilon = \varepsilon$, the latter as defined in Equation~\eqref{eq:epsfam}, under the additional assumption that~\eqref{eq:epscond} holds. This follows from the definition of~$\varepsilon$ together with Lemma~\ref{lem:cControl}, the latter showing that~$0 < \varepsilon (c_1 + c_2) < 1$.
	\end{remark}
	
	\begin{proof}
		Because~$c_1 \in (0, 1)$ by definition, it follows that~$\epsilon c_1 \in (0, \epsilon) \subset (0, 1)$. Furthermore,~$c_2$ is positive, so that~$0 < \epsilon c_2 < 1$ holds (the second inequality is assumed). Therefore, all quantiles appearing in Equations~\eqref{eq:LBLQ}--\eqref{eq:UBUQ} are defined. Due to Lemma~\ref{lem:contdisc}, it is enough to establish the present lemma under the additional assumption that the cdf of~$X_1$ is continuous, \emph{which we shall maintain throughout this proof without further mentioning}.  
		
		We begin by establishing~\eqref{eq:LBLQ}. To this end, let 
		\begin{equation*}
			S_n:=\sum_{i=1}^n\mathds{1}\del[1]{X_i\leq Q_{c_1\epsilon}(X_1)}\qquad \text{and}\qquad \tilde{S}_n:=\sum_{i=1}^n\mathds{1}\del[1]{\tilde{X}_i\leq Q_{c_1\epsilon}(X_1)},
		\end{equation*}
		and note that
		\begin{equation*}
			\cbr[1]{S_n < n(\epsilon-\eta)}
			\subseteq
			\cbr[1]{\tilde{S}_n < n\epsilon }
			\subseteq
			\cbr[1]{\tilde X_{\lceil \epsilon n \rceil}^* 
				\geq
				Q_{c_1\epsilon}(X_1)}.
		\end{equation*}
		Thus, it suffices to show that~$\P\del[1]{S_n\geq n(\epsilon-\eta)}\leq \delta/6$. Noting that~$S_n$ has a Binomial distribution with success probability~$c_1\epsilon\in(0,\epsilon)$, we set up for an application of Part 1.~of Lemma~\ref{lem:Chernoff}.  To this end, note that since~$\epsilon\geq\lambda_1\eta$ and~$c_1<A_+$, it holds that
		\begin{equation*}
			\frac{\epsilon-\eta}{c_1\epsilon}
			=
			\frac{1-\eta/\epsilon}{c_1}
			\geq
			\frac{1-\lambda_1^{-1}\mathds{1}(\eta>0)}{c_1}
			=
			\frac{A_+}{c_1}
			=
			\nu_+(c_1)+1
			>
			1,
		\end{equation*}
		with~$\nu_+$ as defined in Part 1.~of Lemma~\ref{lem:Prophfs}. Therefore, by Part 1.~of Lemma~\ref{lem:Chernoff} 
		\begin{align*}
			\P\del[1]{S_n\geq (\epsilon-\eta)n}
			\leq
			\P\del[1]{S_n\geq (1+\nu_{+}(c_1))c_1\epsilon n}
			\leq 
			e^{-n\epsilon c_1h_+(\nu_{+}(c_1))} &= e^{-n\epsilon f(c_1)} \\ &= \delta/6.
		\end{align*}

		Next, we consider~\eqref{eq:LBUQ}. To this end, redefine 
		\begin{equation*}
			S_n:=\sum_{i=1}^n\mathds{1}\del[1]{X_i\geq Q_{1-c_2\epsilon}(X_1)}\qquad \text{and}\qquad \tilde{S}_n:=\sum_{i=1}^n\mathds{1}\del[1]{\tilde{X}_i\geq Q_{1-c_2\epsilon}(X_1)},
		\end{equation*}
		and note that
		\begin{equation*}
			\cbr[1]{S_n > n(\epsilon+\eta)}
			\subseteq
			\cbr[1]{\tilde{S}_n > n\epsilon }
			\subseteq
			\cbr[1]{\tilde{S}_n \geq \lfloor n\epsilon\rfloor +1 }
			\subseteq
			\cbr[1]{\tilde X_{\lceil(1-\epsilon)n\rceil}^* 
				\geq
				Q_{1-c_2\epsilon}(X_1)};
		\end{equation*}
		the last inclusion using that if at least~$\lfloor \epsilon n \rfloor+1$ of the observations~$\tilde{X}_i$ satisfy~$\tilde{X}_i\geq Q_{1-c_2\epsilon}(X_1)$, then~$\tilde X_{\lceil(1-\epsilon)n\rceil}^*=\tilde{X}_{n-\lfloor \epsilon n\rfloor}^*\geq Q_{1-c_2\epsilon}(X_1)$. Thus, it suffices to show that $\P\del[1]{S_n\leq n(\epsilon+\eta)}\leq \delta/6$. Noting that~$S_n$ has a Binomial distribution with success probability~$c_2\epsilon\in(0,1)$ (it has already been argued that~$c_2\epsilon\in(0,1)$), we set up for an application of Part 2.~of Lemma~\ref{lem:Chernoff}. To this end, note that since~$\epsilon\geq\lambda_1\eta$ and~$c_2>A_-$, it holds that
		\begin{equation*}
			0
			<
			\frac{\epsilon+\eta}{c_2\epsilon}
			=
			\frac{1+\eta/\epsilon}{c_2}
			\leq
			\frac{1+\lambda_1^{-1}\mathds{1}(\eta>0)}{c_2}
			=
			\frac{A_-}{c_2}
			=
			1-\nu_-(c_2)
			<
			1,
		\end{equation*}
		with~$\nu_-$ as defined in Part 2.~of Lemma~\ref{lem:Prophfs}. Therefore, by Part 2.~of Lemma~\ref{lem:Chernoff} 
		\begin{align*}
			\P\del[1]{S_n\leq (\epsilon+\eta)n}
			\leq
			\P\del[1]{S_n\leq (1-\nu_{-}(c_2))c_2\epsilon n}
			\leq 
			e^{-n\epsilon c_2 h_-(\nu_{-}(c_2))} &= e^{-n\varepsilon g(c_2)} \\ &= \delta/6.
		\end{align*}
		
		Next, we consider~\eqref{eq:UBLQ}. To this end, redefine 
		\begin{equation*}
			S_n:=\sum_{i=1}^n\mathds{1}\del[1]{X_i\leq Q_{c_2\epsilon}(X_1)}\qquad \text{and}\qquad \tilde{S}_n:=\sum_{i=1}^n\mathds{1}\del[1]{\tilde{X}_i\leq Q_{c_2\epsilon}(X_1)},
		\end{equation*}
		and note that
		\begin{equation*}
			\cbr[1]{S_n > n(\epsilon+\eta)}
			\subseteq
			\cbr[1]{\tilde{S}_n > n\epsilon }
			\subseteq
			\cbr[1]{\tilde{S}_n \geq \lfloor n\epsilon\rfloor +1 }
			\subseteq
			\cbr[1]{\tilde X_{\lfloor \epsilon n\rfloor+1}^* 
				\leq
				Q_{c_2\epsilon}(X_1)}.
		\end{equation*}
		Thus, it suffices to show that~$\P\del[1]{S_n\leq n(\epsilon+\eta)}\leq \delta/6$. Noting that~$S_n$ has a Binomial distribution with success probability~$c_2\epsilon\in(0,1)$, this has already been established in the proof of the previous case.
		
		Finally, we establish~\eqref{eq:UBUQ}. To this end, redefine 
		\begin{equation*}
			S_n:=\sum_{i=1}^n\mathds{1}\del[1]{X_i\geq Q_{1-c_1\epsilon}(X_1)}\qquad \text{and}\qquad \tilde{S}_n:=\sum_{i=1}^n\mathds{1}\del[1]{\tilde{X}_i\geq Q_{1-c_1\epsilon}(X_1)},
		\end{equation*}
		and note that
		\begin{align*}
			\cbr[1]{S_n < n(\epsilon-\eta)}
			\subseteq
			\cbr[1]{\tilde{S}_n < n\epsilon }
			&\subseteq
			\cbr[1]{\tilde{S}_n \leq \lceil n\epsilon\rceil -1 }\\
			&\subseteq
			\cbr[1]{\tilde X_{\lfloor (1-\epsilon)n\rfloor+1}^* 
				\leq
				Q_{1-c_1\epsilon}(X_1)};
		\end{align*}
		the last inclusion using that if at most~$\lceil \epsilon n\rceil-1$ of the~$\tilde{X}_i$ satisfy that~$\tilde{X}_i\geq Q_{1-c_1\epsilon}(X_1)$ then~$\tilde{X}_{\lfloor (1-\epsilon)n\rfloor+1}^*=\tilde{X}_{n-(\lceil \epsilon n\rceil-1)}^*<Q_{1-c_1\epsilon}(X_1)$. It remains to show that~$\P\del[1]{S_n\geq n(\epsilon-\eta)}\leq \delta/6$. Noting that~$S_n$ has a Binomial distribution with success probability~$c_1\epsilon\in(0,\epsilon)$, this has already been established in the proof of~\eqref{eq:LBLQ}.
	\end{proof}

	\section{Auxiliary results for controlling~$\overline{I}_{n,1}$~$\overline{I}_{n,2}$,~$\overline{I}_{n,3}$ and~$\underline{I}_{n,1}$,~$\underline{I}_{n,2}$,~$\underline{I}_{n,3}$}\label{sec:Is}
	
	The following lemma, which is standard but we could not pinpoint a suitable reference in the literature, bounds the difference between the mean and quantile of a distribution (which is not necessarily continuous). 
	\begin{lemma}\label{lem:quantile_mean}
		Let~$Z$ satisfy~$\sigma_m^m:=\E|Z-\E Z|^m\in[0,\infty)$ for some~$m\in[1,\infty)$. Then, for all~$p\in(0,1)$,
		\begin{equation}\label{eq:QM}
			\E Z-\frac{\sigma_m}{p^{1/m}}\leq Q_p(Z)\leq \E Z+\frac{\sigma_m}{(1-p)^{1/m}}.
		\end{equation}
	\end{lemma}
	
	\begin{proof}
		Fix~$p\in(0,1)$. The statement trivially holds for~$Q_p(Z)=\E Z$, which arises, in particular, if~$\sigma_m = 0$. Thus, let~$Q_p(Z)\neq\E Z$, implying that $\sigma_m\in(0,\infty)$. Denote~$t:=(\E Z-Q_p(Z))/\sigma_m$. 
		
		\emph{Case 1:} If $Q_p(Z)<\E Z$, the second inequality in~\eqref{eq:QM} trivially holds. Elementary properties of the quantile function and Markov's inequality deliver
		\begin{align*}
			p
			\leq
			\P\del[1]{Z\leq Q_p(Z)}
			=
			&\P\del[1]{Z-\E Z\leq Q_p(Z)-\E Z} \\
			&\leq
			\P\del[1]{|Z-\E Z|/\sigma_m\geq |t|} 
			\leq |t|^{-m},
		\end{align*}
		which rearranges to the first inequality in~\eqref{eq:QM}.
		
		\emph{Case 2:}  If~$Q_p(Z)>\E Z$, the first inequality in~\eqref{eq:QM} trivially holds. Elementary properties of the quantile function and Markov's inequality deliver
		\begin{align*}
			1-p
			\leq 
			1-\P\del[1]{Z< Q_p(Z)}
			&=
			\P\del[1]{Z-\E Z\geq Q_p(Z)- \E Z} \\
			&\leq
			\P\del[1]{|Z-\E Z|/\sigma_m \geq |t|} 
			\leq |t|^{-m},
		\end{align*}
		which rearranges to the second inequality in~\eqref{eq:QM}.
	\end{proof}

	In the following we abbreviate~$Q_s=Q_s(X_1)$ for all~$s\in(0,1)$.
	\begin{lemma}\label{lem:noisecontrol}
		Fix~$n\in\N$. Let~$0<s_1<s_2<1$ and Assumption~\ref{ass:setting} be satisfied. Then
		\begin{equation}\label{eq:distcorrupt1}
			\envert[3]{\frac{1}{n}\sum_{i=1}^n\sbr[1]{\phi_{Q_{s_1},Q_{s_2}}(\tilde{X}_i)-\phi_{Q_{s_1},Q_{s_2}}(X_i)}}
			\leq
			\eta \sigma_m\del[3]{\frac{1}{(1-s_2)^{1/m}}+\frac{1}{s_1^{1/m}}}.
		\end{equation}
	\end{lemma}
	\begin{proof}
		Since at most~$\eta n$ observations have been contaminated,
		\begin{align*}
			\envert[3]{\frac{1}{n}\sum_{i=1}^n\sbr[1]{\phi_{Q_{s_1},Q_{s_2}}(\tilde{X}_i)-\phi_{Q_{s_1},Q_{s_2}}(X_i)}}
			&\leq 
			\eta \del[1]{Q_{s_2}-Q_{s_1}} \\
			&\leq	
			\eta \sigma_m\del[3]{\frac{1}{(1-s_2)^{1/m}}+\frac{1}{s_1^{1/m}}} ,
		\end{align*}	
		where the second inequality followed from Lemma~\ref{lem:quantile_mean}.
	\end{proof}
	To establish Lemma~\ref{lem:trimmeanconc} below, we recall Bernstein's inequality from Equation~3.24 of Theorem 3.1.7 in \cite{gine2016mathematical} (note that our statement explicitly requires~$c>0$, which is implicitly imposed in the paragraph preceding their Theorem 3.1.7).
	\begin{theorem}[Bernstein's inequality]\label{thm:bernstein}
		Let~$Z_1, \hdots, Z_n$ be independent centered random variables almost surely bounded by~$c \in (0, \infty)$ in absolute value. Set~$\sigma^2 = n^{-1} \sum_{i = 1}^n \E(Z_i^2)$ and $S_n = \sum_{i = 1}^n Z_i$. Then, $\P(S_n \geq \sqrt{2 n \sigma^2 u} + \frac{cu}{3}) \leq e^{-u}$ for all~$u \geq 0$.
	\end{theorem}
	
	\begin{lemma}\label{lem:trimmeanconc}
		Fix~$n\in\N$ and~$\delta\in(0,1)$. Let~$0<s_1<s_2<1$ and Assumption~\ref{ass:setting} be satisfied. Let
		\begin{equation*}
			\tau:=\del[3]{\frac{\sigma_m}{(1-s_2)^{1/m}}+\frac{\sigma_m}{s_1^{1/m}}}^{2-(m\wedge 2)}\sigma_{m\wedge 2}^{m\wedge 2}.
		\end{equation*}
		Then each of 
		\begin{equation}
		\begin{aligned}\label{eq:Bernstein1}
			&\frac{1}{n}\sum_{i=1}^n\sbr[1]{\phi_{Q_{s_1},Q_{s_2}}(X_i)-\E\phi_{Q_{s_1},Q_{s_2}}(X_i)} \\
			&\hspace{2cm} \geq 
			-\sqrt{\frac{2\tau\log(6/\delta)}{n}}-  \del[3]{\frac{\sigma_m}{(1-s_2)^{1/m}}+\frac{\sigma_m}{s_1^{1/m}}}\frac{\log(6/\delta)}{3n}
		\end{aligned}
		\end{equation}
		and
		\begin{equation}
		\begin{aligned}\label{eq:Bernstein2}
			&\frac{1}{n}\sum_{i=1}^n\sbr[1]{\phi_{Q_{s_1},Q_{s_2}}(X_i)-\E\phi_{Q_{s_1},Q_{s_2}}(X_i)} \\
			&\hspace{2cm} \leq 
			\sqrt{\frac{2\tau\log(6/\delta)}{n}}+  \del[3]{\frac{\sigma_m}{(1-s_2)^{1/m}}+\frac{\sigma_m}{s_1^{1/m}}}\frac{\log(6/\delta)}{3n}
		\end{aligned}
		\end{equation}
		holds with probability at least~$1-\delta/6$.
	\end{lemma}
	\begin{proof}
		The statement is trivially true in case~$\sigma_m = 0$ (which implies~$Q_{s_1} = Q_{s_2}$). Hence, we shall assume throughout that~$\sigma_m > 0$. We first make two observations that will allow us to apply Bernstein's inequality. For~$i=1,\hdots, n$, note that
		\begin{align*}
			Y_i&:=\envert[1]{\phi_{Q_{s_1},Q_{s_2}}(X_i)-\E\phi_{Q_{s_1},Q_{s_2}}(X_i)} \\ &\leq  
			Q_{s_2}-Q_{s_1}
			\leq 
			\del[3]{\frac{\sigma_m}{(1-s_2)^{1/m}}+\frac{\sigma_m}{s_1^{1/m}}}\in(0,\infty),
		\end{align*}
		where the second inequality followed from Lemma~\ref{lem:quantile_mean}. Therefore,
		\begin{equation*}
			\E Y_1^2
			=
			\E\del[1]{|Y_1|^{2-(m\wedge 2)}|Y_1|^{m\wedge 2}}
			\leq
			\del[3]{\frac{\sigma_m}{(1-s_2)^{1/m}}+\frac{\sigma_m}{s_1^{1/m}}}^{2-(m\wedge 2)}\E|Y_1|^{m\wedge 2}
			\leq \tau,
		\end{equation*}
		where the last inequality used that~$\E|Y_1|^k \leq \E|X_1-\mu|^k = \sigma_k^k$ for~$k=m\wedge 2$, cf., e.g., Corollary 3 in~\cite{cstudd}. 
		
		Now, standard arguments combined with Bernstein's inequality (Theorem~\ref{thm:bernstein}) show that~\eqref{eq:Bernstein1} and~\eqref{eq:Bernstein2}, respectively, holds with probability at least~$1-\delta/6$.
	\end{proof}
	
	\begin{lemma}\label{lem:meancontrol}
		Let~$0<s_1<s_2<1$ and Assumption~\ref{ass:setting} be satisfied. Then
		\begin{equation}\label{eq:mean1}
			\E\phi_{Q_{s_1},Q_{s_2}}(X_1)-\mu
			\geq
			-2\sigma_ms_1^{1-\frac{1}{m}}-\sigma_m\del[2]{1+\sbr[2]{\frac{1-s_2}{s_2}}^{\frac{1}{m}}}(1-s_2)^{1-\frac{1}{m}},
		\end{equation}	
		and 
		\begin{equation}\label{eq:mean2}
			\E\phi_{Q_{s_1},Q_{s_2}}(X_1)-\mu
			\leq
			2\sigma_m(1-s_2)^{1-\frac{1}{m}}+\sigma_m\del[2]{1+\sbr[2]{\frac{s_1}{1-s_1}}^{\frac{1}{m}}}s_1^{1-\frac{1}{m}}.
		\end{equation}
	\end{lemma}
	\begin{proof}
		We write~$\phi_{Q_{s_1},Q_{s_2}}(X_1)-\mu$ equivalently as
		\begin{equation*}
			(X_1-\mu)\mathds{1}(Q_{s_1}\leq X_1\leq Q_{s_2})
			+(Q_{s_1}-\mu)\mathds{1}(X_1<Q_{s_1})
			+(Q_{s_2}-\mu)\mathds{1}(Q_{s_2}<X_1),
		\end{equation*}
		such that~$\E\phi_{Q_{s_1},Q_{s_2}}(X_1)-\mu$ equals
		\begin{align}
			&\E\del[1]{(X_1-\mu)\mathds{1}(Q_{s_1}\leq X_1\leq Q_{s_2})}
			+
			(Q_{s_1}-\mu)\P(X_1<Q_{s_1})
			+
			(Q_{s_2}-\mu)\P(X_1>Q_{s_2})\notag\\
			&=  -\E(X_1-\mu)\mathds{1}(X_1<Q_{s_1})
			-\E(X_1-\mu)\mathds{1}(X_1>Q_{s_2})
			+
			(Q_{s_1}-\mu)\P(X_1<Q_{s_1})\notag\\
			&+
			(Q_{s_2}-\mu)\P(X_1>Q_{s_2})
			\label{eq:identical}.
		\end{align}
		We now establish~\eqref{eq:mean1}. Using H{\"o}lder's inequality (with the usual conventions in case~$m=1$) to bound the first two summands on the right-hand side of~\eqref{eq:identical}, and Lemma~\ref{lem:quantile_mean} to bound the last two summands, along with~$\P(X_1<Q_{s_1})\leq s_1$ and~$\P(X_1>Q_{s_2})=1-\P(X_1\leq Q_{s_2})\leq 1-s_2$, it follows that
		\begin{align*}
			\E\phi_{Q_{s_1},Q_{s_2}}(X_1)-\mu
			&\geq 
			-\sigma_m s_1^{1-\frac{1}{m}}-\sigma_m(1-s_2)^{1-\frac{1}{m}}-\frac{\sigma_m}{s_1^{1/m}}s_1-\frac{\sigma_m}{s_2^{1/m}}(1-s_2)\\
			&= 
			-2\sigma_ms_1^{1-\frac{1}{m}}-\sigma_m\del[2]{1+\sbr[2]{\frac{1-s_2}{s_2}}^{\frac{1}{m}}}(1-s_2)^{1-\frac{1}{m}}.
		\end{align*}
		To prove~\eqref{eq:mean2}, we use~\eqref{eq:identical} and the same inequalities as above to conclude that
		\begin{align*}
			\E\phi_{Q_{s_1},Q_{s_2}}(X_1)-\mu
			&\leq 
			\sigma_m s_1^{1-\frac{1}{m}}+\sigma_m(1-s_2)^{1-\frac{1}{m}}+\frac{\sigma_m}{(1-s_1)^{\frac{1}{m}}}s_1+\frac{\sigma_m}{(1-s_2)^{\frac{1}{m}-1}}\\
			&= 
			2\sigma_m(1-s_2)^{1-\frac{1}{m}}+\sigma_m\del[2]{1+\sbr[2]{\frac{s_1}{1-s_1}}^{\frac{1}{m}}}s_1^{1-\frac{1}{m}}.
		\end{align*}
	\end{proof}

	\section{Proof of Theorem~\ref{thm:maintext}}
	
	Recall that throughout~$c_1$ and~$c_2$ are as defined in~\eqref{eq:c1} and~\eqref{eq:c2}, respectively. By Lemma~\ref{lem:quantiles} together with Remark~\ref{rem:qlepschoice} and the arguments leading up to \eqref{eq:lbub}--\eqref{eq:decomp2}, one has with probability at least~$1-\frac{4}{6}\delta$ that 
	\begin{equation*}
		\envert[1]{\hat{\mu}_{n}(\eps(\eta))-\mu}
		\leq
		\del[1]{\overline I_{n,1}+\overline I_{n,2}+\overline I_{n,3}} \vee -\del[1]{\underline I_{n,1}+\underline I_{n,2}+\underline I_{n,3}}.
	\end{equation*}	
	In the following, we employ Lemmas~\ref{lem:noisecontrol}, \ref{lem:trimmeanconc}, and~\ref{lem:meancontrol}, with~$s_1=c_2\eps$ and~$s_2=1-c_1\eps$, to bound~$\overline I_{n,1}+\overline I_{n,2}+\overline I_{n,3}$ from above.\footnote{Identical arguments based on~$s_1=c_1\eps$ and~$s_2=1-c_2\eps$ establish the same upper bounds on~$-\underline{I}_{n,i}$ instead of~$\overline{I}_{n,i}$, respectively, for~$i = 1, 2, 3$. We omit the details.\label{foot:lowerbar}} By~\eqref{eq:epscond} and Lemma~\ref{lem:cControl} it follows that~$s_1<s_2$ as required in these lemmas. We define, for positive real numbers~$d_1$ and~$d_2$, 
	\begin{equation*}
		A_m(d_1,d_2):=\frac{1}{d_1^{1/m}}+\frac{1}{d_2^{1/m}}\quad\text{and}\quad B_m(d_1,d_2):=2d_1^{1-\frac{1}{m}}+\sbr[2]{1+\del[2]{\frac{d_2}{d_1}}^{\frac{1}{m}}}d_2^{1-\frac{1}{m}},
	\end{equation*}
	
	If~$\eta=0$, then~$\overline{I}_{n,1}=0$ as well. If~$\eta \neq 0$, then, by Lemma~\ref{lem:noisecontrol} and~$\eps \geq \lambda_1\eta$, we have
	\begin{equation*}
		\overline{I}_{n,1}
		\leq
		\eta \sigma_m\del[3]{\frac{1}{(c_1\eps)^{1/m}}+\frac{1}{(c_2\eps)^{1/m}}} \leq \sigma_m \lambda_1^{-\frac{1}{m}} A_m(c_1, c_2)   \eta^{1-\frac{1}{m}}.
	\end{equation*}
	
	Next, by Lemma~\ref{lem:trimmeanconc} and~$\eps\geq \lambda_2\frac{\log(6/\delta)}{n}$, it holds with probability at least~$1-\delta/6$ (the ``final''~$1-\delta/6$ comes from bounding~$-\underline{I}_{n,2}$ by identical arguments, cf.~Footnote~\ref{foot:lowerbar}) that in case~$m\geq 2$ (where~$\tau$ in Lemma~\ref{lem:trimmeanconc} equals~$\sigma_2^2$):
	\begin{align*}
		\overline I_{n,2}
		&\leq
		\sqrt{\frac{2\sigma_2^2\log(6/\delta)}{n}}+  \del[3]{\frac{\sigma_m}{(c_1\eps)^{1/m}}+\frac{\sigma_m}{(c_2\eps)^{1/m}}}\frac{\log(6/\delta)}{3n}\\
		&\leq \sqrt{2}\sigma_m
		\sqrt{\frac{\log(6/\delta)}{n}}+ \sigma_m \lambda_2^{-\frac{1}{m}} (A_m(c_1, c_2)/3) \sbr[2]{\frac{\log(6/\delta)}{n}}^{1-\frac{1}{m}} \\
		&\leq \sigma_m  \cdot \left\{ \sqrt{2} + \lambda_2^{-\frac{1}{m}} A_m(c_1, c_2)/3 \right\}   \cdot	\sqrt{\frac{\log(6/\delta)}{n}},
	\end{align*}
	the last inequality following from~$\log(6/\delta)/n<1$ by~\eqref{eq:epscond}. In the case where~$m\in[1,2)$, the quantity~$\tau$ in Lemma~\ref{lem:trimmeanconc} equals~$\sigma_m^2\del[2]{\frac{1}{(1-s_2)^{1/m}}+\frac{1}{s_1^{1/m}}}^{2-m}$, such that with probability at least~$1-\delta/6$, using similar arguments as in the previous case, particularly~$\eps\geq \lambda_2\frac{\log(6/\delta)}{n}$,
	\begin{align*}
		\overline I_{n,2}
		\leq
		\sigma_m \left\{
		\sqrt{2}  \del[3]{\lambda_2^{-\frac{1}{m}} A_m(c_1, c_2)}^{1-\frac{m}{2}}+ \lambda_2^{-\frac{1}{m}} \frac{A_m(c_1, c_2)}{3}  \right\}   \cdot  \left[\frac{\log(6/\delta)}{n}\right]^{1-\frac{1}{m}}.
	\end{align*}
	We can summarize both cases in the following way
	\begin{equation*}
		\overline I_{n,2} \leq 	\sigma_m \left\{
		\sqrt{2}  \del[3]{\lambda_2^{-\frac{1}{m}} A_m(c_1, c_2)}^{1-\frac{m \wedge 2}{2}}+ \lambda_2^{-\frac{1}{m}} \frac{A_m(c_1, c_2)}{3}  \right\}  \cdot \left[\frac{\log(6/\delta)}{n}\right]^{1-\frac{1}{m \wedge 2}}.
	\end{equation*}
	Finally, by Lemma~\ref{lem:meancontrol}, and using that by Lemma~\ref{lem:cControl} and~\eqref{eq:epscond} it holds that $\eps(c_1+c_2)<1$ such that~$1-c_2\eps >c_1\eps$, we obtain
	\begin{align*}
		\overline I_{n,3}
		&\leq
		2\sigma_m(c_1\eps)^{1-\frac{1}{m}}+\sigma_m\del[2]{1+\sbr[2]{\frac{c_2\eps}{1-c_2\eps}}^{\frac{1}{m}}}(c_2\eps)^{1-\frac{1}{m}}\\
		&\leq 
		2\sigma_m(c_1\eps)^{1-\frac{1}{m}}+\sigma_m\del[2]{1+\sbr[2]{\frac{c_2\eps}{c_1\eps}}^{\frac{1}{m}}}(c_2\eps)^{1-\frac{1}{m}}\\
		&= \sigma_m \eps^{1-\frac{1}{m}}  \cdot B_m(c_1, c_2)  \\
		&\leq 
		\sigma_m  B_m(c_1, c_2)  \cdot \left[ \lambda_1 ^{1-\frac{1}{m}}\cdot \eta^{1-\frac{1}{m}} +\lambda_2^{1-\frac{1}{m}}\cdot \left[\frac{\log(6/\delta)}{n}\right]^{1-\frac{1}{m\wedge 2}}  \right],
	\end{align*}
	the last inequality using sub-additivity of~$z\mapsto z^{1-\frac{1}{m}}$ (recalling again that $\log(6/\delta)/n<1$).
	
	Summarizing (cf.~also Footnote~\ref{foot:lowerbar}), with probability at least~$1-\delta$ we obtain the following upper bound on~$\del[1]{\overline I_{n,1}+\overline I_{n,2}+\overline I_{n,3}} \vee -\del[1]{\underline I_{n,1}+\underline I_{n,2}+\underline I_{n,3}}$ (and hence on~$\envert[1]{\hat{\mu}_{n}(\eps(\eta))-\mu}$):
	\begin{align*}
		&\sigma_m \lambda_1^{-\frac{1}{m}} A_m(c_1, c_2)  \eta^{1-\frac{1}{m}} \\
		& +  \sigma_m \left\{
		\sqrt{2}  \del[3]{\lambda_2^{-\frac{1}{m}} A_m(c_1, c_2)}^{1-\frac{m \wedge 2}{2}}+ \lambda_2^{-\frac{1}{m}} A_m(c_1, c_2) /3  \right\}   \cdot \left[\frac{\log(6/\delta)}{n}\right]^{1-\frac{1}{m \wedge 2}} \\
		& +  \sigma_m  B_m(c_1, c_2)  \cdot \left[ \lambda_1 ^{1-\frac{1}{m}}\cdot \eta^{1-\frac{1}{m}} +\lambda_2^{1-\frac{1}{m}}\cdot \left[\frac{\log(6/\delta)}{n}\right]^{1-\frac{1}{m\wedge 2}}  \right], 
	\end{align*}
	which, collecting terms, re-arranges to
	\begin{equation*}
		\sigma_m  \cdot \left[ \mathfrak{A}^{\dagger}_m(c_1, c_2) \cdot \eta^{1-\frac{1}{m}} + \mathfrak{B}^{\dagger}_m(c_1, c_2)  \cdot \left[\frac{\log(6/\delta)}{n}\right]^{1-\frac{1}{m\wedge 2}} \right],
	\end{equation*}
	with
	\begin{equation*}
		\mathfrak{A}^{\dagger}_m(c_1, c_2) := \lambda_1^{-\frac{1}{m}} \cdot \left[ A_m(c_1, c_2) + \lambda_1 B_m(c_1, c_2) \right],
	\end{equation*}
	and
	\begin{align*}
		\mathfrak{B}^{\dagger}_m(c_1, c_2) &:= 
		\sqrt{2}  \del[3]{\lambda_2^{-\frac{1}{m}} A_m(c_1, c_2)}^{1-\frac{m \wedge 2}{2}}\\
		&+ \lambda_2^{-\frac{1}{m}} \left( (A_m(c_1, c_2)/3)    + \lambda_2 B_m(c_1, c_2) \right).
	\end{align*}
	Recall from Lemma~\ref{lem:cControl} the following bounds
	\begin{align*}
		& \mathfrak{l}(\lambda_1, \lambda_2) := (1-\lambda_1^{-1})\exp \del[2]{{-\frac{1}{\lambda_2(1-\lambda_1^{-1})}-1}} \leq  c_1\leq 1,\\
		&1 \leq c_2 \leq 2+\lambda_2^{-1}+\sqrt{\lambda_2^{-2}+4\lambda_2^{-1}} =: \mathfrak{u}(\lambda_1, \lambda_2).
	\end{align*}
	It hence follows that
	\begin{equation*}
		A_m(c_1, c_2) \leq A_m(\mathfrak{l}(\lambda_1, \lambda_2), 1)
	\end{equation*}
	and that
	\begin{equation*}
		B_m(c_1, c_2) 
		\leq 
		2+\sbr[2]{1+\del[2]{\frac{\mathfrak{u}(\lambda_1, \lambda_2)}{\mathfrak{l}(\lambda_1, \lambda_2)}}^{\frac{1}{m}}}\mathfrak{u}(\lambda_1, \lambda_2)^{1-\frac{1}{m}} =: \overline{B}_m(\lambda_1, \lambda_2), 
	\end{equation*}
	from which we can conclude that with probability at least~$1-\delta$, it holds that
	\begin{equation*}
		\envert[1]{\hat{\mu}_{n}(\eps(\eta))-\mu} \leq \sigma_m  \cdot \left[ \mathfrak{A}_m(\lambda_1, \lambda_2) \cdot \eta^{1-\frac{1}{m}} + \mathfrak{B}_m(\lambda_1, \lambda_2) \cdot \left[\frac{\log(6/\delta)}{n}\right]^{1-\frac{1}{m\wedge 2}} \right],
	\end{equation*}
	where
	\begin{align*}
		\mathfrak{A}_m(\lambda_1, \lambda_2) &:= \lambda_1^{-\frac{1}{m}} \cdot \left[ A_m(\mathfrak{l}(\lambda_1, \lambda_2), 1) + \lambda_1 \overline{B}_m(\lambda_1, \lambda_2) \right]\\
		\mathfrak{B}_m(\lambda_1, \lambda_2) &:= \sqrt{2}  \del[3]{\lambda_2^{-\frac{1}{m}} A_m(\mathfrak{l}(\lambda_1, \lambda_2), 1)}^{1-\frac{m \wedge 2}{2}}\\
		&+ \lambda_2^{-\frac{1}{m}} \left( (A_m(\mathfrak{l}(\lambda_1, \lambda_2), 1)/3)   + \lambda_2 \overline{B}_m(\lambda_1, \lambda_2) \right). 
	\end{align*}
	
	The statement in Footnote~\ref{foot:m1eta0} follows from a simple adaptation of the above argument to the case~$m = 1$ and~$\eta = 0$.
	
	\section{Proof of Theorem~\ref{thm:adapt}}
	\begin{proof}[Proof of Theorem~\ref{thm:adapt}]
		We first argue that~$\hat{\mu}_{n,A}$ is well-defined. By assumption,~$\eps_A(\eta_{g^*})$ satisfies~\eqref{eq:epscondLepski}, such that~$\mathbb{I}(\eta_{g^*})=\mathbb{B}\del[1]{\hat{\mu}_n(\eps_A(\eta_{g^*})),B(\eta_{g^*})}$. Thus, on the one hand, if~$\hat{g}=g_{\max}$, then~$\bigcap_{j=1}^{\hat g} \mathbb{I}(\eta_j)$ is a non-empty finite interval [as it intersects over the finite interval $\mathbb{B}\del[1]{\hat{\mu}_n(\eps_A(\eta_{g^*})),B(\eta_{g^*})}$]. If, on the other hand,~$\hat{g}<g_{\max}$, then~$\bigcap_{j=1}^{\hat g+1}\mathbb{I}(\eta_j)=\emptyset$ by definition of~$\hat{g}$. Thus,~$\bigcap_{j=1}^{\hat g}\mathbb{I}(\eta_j)\neq \R$, and it follows that~$\mathbb{I}(\eta_j)=\mathbb{B}\del[1]{\hat{\mu}_n(\eps_A(\eta_j)),B(\eta_j)}$ for at least one~$j=1,\hdots,\hat{g}$. Thus, $\bigcap_{j=1}^{\hat g}\mathbb{I}(\eta_j)$ is again a non-empty finite interval, and its midpoint~$\hat{\mu}_{n,A}$ is well-defined.
		
		We now establish~\eqref{eq:UBadapt}. Let~$j\in[g^*]=\cbr[1]{1,\hdots,g^*}$, such that~$\eta_{\min}\leq \eta_j$. If, in addition,~$\eps_A(\eta_j)$ satisfies~\eqref{eq:epscondLepski}, then~$\mathbb{I}(\eta_j)=\mathbb{B}\del[1]{\hat{\mu}_n(\eps_A(\eta_j)),B(\eta_j)}$, and it holds by Theorem~\ref{thm:maintext} that $\mu\in\mathbb{I}(\eta_j)$ with probability at least~$1-\delta/g_{\max}$. If~$\eps_A(\eta_j)$ does not satisfy~\eqref{eq:epscondLepski} then~$\mathbb{I}(\eta_j)=\R$ and~$\mu\in \mathbb{I}(\eta_j)$ with probability one. Thus, by the union bound,
		\begin{equation*}
			\mu \in \bigcap_{j=1}^{g^*}\mathbb{I}(\eta_j) \qquad\text{with probability at least }1-\delta.
		\end{equation*}
		On~$\cbr[1]{\mu \in \bigcap_{j=1}^{g^*}\mathbb{I}(\eta_j)}$, which we shall suppose to occur in what follows, it holds that~$\hat{g}\geq g^*$, such that also
		\begin{equation*}
			\hat{\mu}_{n,A}\in \bigcap_{j=1}^{\hat g} \mathbb{I}(\eta_j)
			\subseteq 
			\bigcap_{j=1}^{g^*}\mathbb{I}(\eta_j). 
		\end{equation*}
		Thus,~$\hat{\mu}_{n,A}$ and~$\mu$ both belong to
		\begin{equation*}
			\bigcap_{j=1}^{g^*}\mathbb{I}(\eta_j)	
			\subseteq
			\mathbb{I}(\eta_{g^*})
			=
			\mathbb{B}\del[1]{\hat{\mu}_n(\eps_A(\eta_{g^*})),B(\eta_{g^*})},
		\end{equation*}
		where we used that~$\eps_A(\eta_{g^*})$ satisfies~\eqref{eq:epscondLepski}. It follows that
		\begin{equation}\label{eq:forremark}
			|\hat{\mu}_{n,A}-\mu|
			\leq 2B(\eta_{g^*}).
		\end{equation}
		
		In case~$g^* < g_{\max}$, it holds that~$\rho\eta_{g^*} < \eta_{\min} \leq \eta_{g^*}$. Since~$z\mapsto B(z)$ is non-decreasing,~$B(\eta_{g^*})$ is then bounded from above by~$B\del[1]{\frac{\eta_{\min}}{\rho}}$, which equals
		\begin{align*}
			\sigma_m \cdot \left(
			\mathfrak{A}_m(\lambda_1, \lambda_2) \cdot \left[\frac{\eta_{\min}}{\rho}\right]^{1-\frac{1}{m}}+\mathfrak{B}_m(\lambda_1, \lambda_2) \cdot \del[2]{\frac{\log(6g_{\max}/\delta)}{n}}^{1-\frac{1}{m\wedge 2}}\right).
		\end{align*}
		
		In case~$g^*=g_{\max}=\lceil \log_\rho(2\log(6/\delta)/n)\rceil$, it follows that
		\begin{equation*}
			\eta_{g^*} = \eta_{g_{\max}} = 0.5 \rho^{g_{\max}} \leq \log(6/\delta)/n \leq \del[2]{\frac{\log(6g_{\max}/\delta)}{n}},
		\end{equation*}
		and we recall that~$\log(6g_{\max}/\delta)/n < 1$ as a consequence of the assumption that~$\eps_A(\eta_{g^*})$ satisfies~\eqref{eq:epscondLepski}. Thus, in this case~$B(\eta_{g^*})$ is bounded from above by
		\begin{align*}
			 \sigma_m \cdot \left(
			\mathfrak{A}_m(\lambda_1, \lambda_2) \cdot \del[2]{\frac{\log(6g_{\max}/\delta)}{n}}^{1-\frac{1}{m}}+\mathfrak{B}_m(\lambda_1, \lambda_2) \cdot \del[2]{\frac{\log(6g_{\max}/\delta)}{n}}^{1-\frac{1}{m\wedge 2}}\right) \\
			 \leq \sigma_m \cdot \left(
			\mathfrak{A}_m(\lambda_1, \lambda_2)  +\mathfrak{B}_m(\lambda_1, \lambda_2) \right) \cdot  \del[2]{\frac{\log(6g_{\max}/\delta)}{n}}^{1-\frac{1}{m\wedge 2}}.
		\end{align*}
		Combining the two cases, we obtain the claimed bound.
	\end{proof}
	\begin{remark}\label{rem:altadapt}
		The alternative estimator~$\tilde{\mu}_{n}=\hat{\mu}_n(\eps_A(\eta_{\hat{g}}))$ in Remark~\ref{rem:otherestim} obeys the following performance guarantee. As argued in the proof of Theorem~\ref{thm:adapt} above (with all notation as there),
		\begin{equation*}
			\mu \in \bigcap_{j=1}^{g^*}\mathbb{I}(\eta_j) \qquad\text{with probability at least }1-\delta.
		\end{equation*}and on this event~$\hat{g}\geq g^*$. Thus, 
		\begin{equation*}
			\emptyset\neq\bigcap_{j=1}^{\hat g} \mathbb{I}(\eta_j)
			\subseteq 
			\bigcap_{j=1}^{g^*}\mathbb{I}(\eta_j). 
		\end{equation*}
		Next,~$\eps_A(\eta_{\hat g})\leq \eps_A(\eta_{g^*})$ with~$\eps_A(\eta_{g^*})$ and hence~$\eps_A(\eta_{\hat g})$ satisfying~\eqref{eq:epscondLepski} (the former by assumption) such that~$\mathbb{I}(\eta_{\hat g})=\mathbb{B}\del[1]{\hat{\mu}_n(\eps_A(\eta_{\hat g})),B (\eta_{\hat g})}$ and~$\mathbb{I}(\eta_{g^*})=\mathbb{B}\del[1]{\hat{\mu}_n(\eps_A(\eta_{g^*})),B (\eta_{g^*})}$. Thus, denoting by~$\hat{y}$ an element of the left intersection in the previous display, it holds that~$\hat y\in \mathbb{B}\del[1]{\hat{\mu}_n(\eps_A(\eta_{\hat g})),B (\eta_{\hat g})}$ and~$\hat y\in \mathbb{B}\del[1]{\hat{\mu}_n(\eps_A(\eta_{g^*})),B (\eta_{g^*})}$. By the triangle inequality~$\tilde{\mu}_{n}=\hat{\mu}_n(\eps_A(\eta_{\hat{g}}))$ hence satisfies
		\begin{align}\label{eq:closetooracle}
			\envert[1]{\tilde{\mu}_{n}-\hat{\mu}_n(\eps_A(\eta_{g^*}))}
			&\leq
			|\hat{\mu}_n(\eps_A(\eta_{\hat{g}}))-\hat{y}|+|\hat{y}-\hat{\mu}_n(\eps_A(\eta_{g^*}))|
			\leq
			B(\eta_{\hat{g}})+B(\eta_{g^*})\notag\\
			&\leq
			2B(\eta_{g^*}).
		\end{align}
		In addition, since~$\mu\in\mathbb{I}(\eta_{g^*})=\mathbb{B}\del[1]{\hat{\mu}_n(\eps_A(\eta_{g^*})),B (\eta_{g^*})}$ we have $|\hat{\mu}_n(\eps_A(\eta_{g^*}))-\mu|\leq B(\eta_{g^*})$. In combination with the previous display, this yields~$|\tilde{\mu}_{n}-\mu|\leq 3B(\eta_{g^*})$. Splitting into the cases of~$g^*<g_{\max}$ and $g^*=g_{\max}$ like at the end of the proof of Theorem~\ref{thm:adapt}, we conclude as in the arguments commencing from~\eqref{eq:forremark}.
	\end{remark}
	
\end{appendix}

\end{document}